\documentclass{article}

\usepackage{mathrsfs}

\usepackage{amsmath}

\usepackage{amssymb}
\usepackage{enumitem}
\usepackage[dvips]{graphicx}

\usepackage{amsthm}

\allowdisplaybreaks[4]

\makeatletter

 \@addtoreset{equation}{section}
\makeatother

\newtheorem{dfn}{Definition}[section]
\newtheorem{thm}[dfn]{Theorem}
\newtheorem{prop}[dfn]{Proposition}
\newtheorem{lem}[dfn]{Lemma}
\newtheorem{rem}[dfn]{Remark}

\newcommand{\no}[3]{\| #1\|^{#2}_{#3}}	
\newcommand{\inn}[2]{\langle\!\langle #1,#2\rangle\!\rangle}		
\newcommand{\mean}[1]{\langle\!\langle #1\rangle\!\rangle}		
\newcommand{\trans}{{}^{\top}}	
\newcommand{\nom}[3]{[\![ #1 ]\!]^{#2}_{#3}}

\begin{document}
\title{\large
\bf{On the spectral properties for the linearized problem around space-time periodic states of the compressible Navier-Stokes equations}}
\author{Mohamad Nor Azlan${}^1$, Shota Enomoto${}^2$, Yoshiyuki Kagei${}^3$}
\date{}

\footnotetext[1]{Kuala Lumpur, MALAYSIA (mohamadnorazlan189@gmail.com)}
\footnotetext[2]{General Education Department, National Institute of Technology, Toba College, Mie 517-8501, JAPAN (enomoto-s@toba-cmt.ac.jp)}
\footnotetext[3]{Department of Mathematics, Tokyo Institute of Technology, Tokyo 152-8551, JAPAN (kagei@math.titech.ac.jp)}

\maketitle

\begin{abstract}
This paper studies the linearized problem for the compressible Navier-Stokes equation around space-time periodic state in an infinite layer of $\mathbb{R}^n$ ($n=2,3$), and the spectral properties of the linearized evolution operator is investigated. 
It is shown that if the Reynolds and Mach numbers are sufficiently small, 
then the asymptotic expansions of the Floquet exponents near the imaginary axis for the Bloch transformed linearized problem are obtained for small Bloch parameters, which would give the asymptotic leading part of the linearized solution operator as $t\rightarrow\infty$. 
\vspace{0.3cm}

\end{abstract}

\section{Introduction}
\numberwithin{equation}{section}
This paper is concerned with the stability of space-time periodic state of the system of equations for a barotropic motion of a viscous and compressible fluid 
\begin{gather}
\partial_{\tilde{t}}\tilde{\rho}+\mathrm{div}_{\tilde{x}}(\tilde{\rho} \tilde{v})=0, \label{den1} \\
\tilde{\rho}(\partial_{\tilde{t}}\tilde{v}+\tilde{v}\cdot\nabla_{\tilde{x}} \tilde{v})-\mu\Delta_{\tilde{x}} \tilde{v}-(\mu+\mu')\nabla_{\tilde{x}}\mathrm{div}_{\tilde{x}}\tilde{v}+\nabla_{\tilde{x}} \tilde{p}(\tilde{\rho})=\tilde{\rho} \tilde{G} \label{vel1}
\end{gather}
in an $n$ dimensional infinite layer $\tilde{\Omega}=\mathbb{R}^{n-1} \times (0,d),\text{ for } n=2,3$.
Here, $\tilde{\rho}=\tilde{\rho}(\tilde{x},\tilde{t})$ and $\tilde{v}=\trans(\tilde{v}_1(\tilde{x},\tilde{t}),\cdots,\tilde{v}_n(\tilde{x},\tilde{t}))$ denote the unknown density and the velocity at time $\tilde{t}\geq 0$ and position $\tilde{x} \in \tilde{\Omega}$, respectively. $\tilde{p}(\tilde{\rho})$ is the pressure; we assume that $\tilde{p}(\tilde{\rho})$ is a smooth function of $\tilde \rho$ and satisfies 
$$\tilde{p}'(\rho_*)=\frac{d \tilde{p}}{d \tilde{\rho}}(\rho_*)>0$$
 for a given constant $\rho_* > 0$. $\mu$ and $\mu'$ are the viscosity coefficients; we assume that $\mu$ and $\mu'$ are constants and the shear viscosity $\mu$ is positive and the bulk viscosity $\frac{2}{n}\mu + \mu'$ is nonnegative. The system \eqref{den1}-\eqref{vel1} is classified in quasilinear hyperbolic-parabolic systems. 

Due to a technical reason, we also assume that $\frac{\mu'}{\mu}$ satisfies 
\begin{equation}
\frac{\mu'}{\mu}\leq \mu_* \label{assump mu}
\end{equation}
for a given constant $\mu_*>0$. (See Remark \ref{upperbound}.) 
$\tilde{G}= \tilde{G}(\tilde{x},\tilde{t})$ is a given external force satisfying
\begin{equation}\label{force}
\tilde{G}(\tilde{x}'+\frac{2\pi}{\tilde{\alpha}_i}\textbf{e}_i',\tilde{x}_n,\tilde{t})=\tilde{G}(\tilde{x}',\tilde{x}_n,\tilde{t}), \ \tilde{G}(\tilde{x}',\tilde{x}_n,\tilde{t}+T)=\tilde{G}(\tilde{x}',\tilde{x}_n,\tilde{t})
\end{equation}
for all $\tilde x'\in \mathbb{R}^{n -1}$, $\tilde x_{n} \in (0, d)$ and $\tilde t \in \mathbb{R}$, where $\tilde{\alpha}_i ~(i=1,\cdots,n-1)$ are positive constants and $\mathbf{e}_i'=\trans(0,\cdots,0,\stackrel{i}{1},0,\cdots,0)\in\mathbb{R}^{n-1}$.
The system \eqref{den1}-\eqref{vel1} is considered under the boundary condition and initial condition  
\begin{gather} \label{cond}
\tilde{v}|_{\tilde{x}_n =0,d}=0 , \\ \label{cond1}
(\tilde{\rho},\tilde{v})|_{\tilde{t}=0}= (\tilde{\rho}_0,\tilde{v}_0). 
\end{gather}

One can see that if $\tilde{G}$ is sufficiently small, the system (\ref{den1})-(\ref{vel1}) with (\ref{cond}) has a space-time periodic state $\tilde{u}_p=\trans(\tilde{\rho}_p, \tilde{v}_p)$.
The purpose of this paper is to investigate the spectral properties of the linearized evolution operator around the space-time periodic state $\tilde{u}_p$ which will be useful in the study of the large time behavior of solutions around $\tilde{u}_{p}$. 

If the external force $\tilde G$ takes the form $\tilde G = \tilde G_{para} = \trans(\tilde g^{1}(\tilde x_{n}, \tilde t), 0, \cdots, 0)$, then the system (1.1)-(1.2) with the boundary condition \eqref{cond} has a time-periodic parallel flow, i.e., a time periodic solution of the form $\tilde u_{para}= \trans(\rho_{*}, \tilde v_{para})$ with $\tilde v_{para}=\trans(\tilde v_{para}^{1}(\tilde x_{n}, \tilde t), 0, \cdots,0)$. The stability of parallel flows has been widely studied in the hydrodynamic stability theory. 
As for the mathematical study of the stability of time periodic parallel flows of \eqref{den1}-\eqref{vel1}, the nonlinear dynamics of solutions  
around time periodic parallel flows was investigated by Brezina \cite{Brezina3}. 
(See also \cite{Brezina2, Brezina1} for the linearized analysis.)
It was proved in \cite{Brezina3} that if the Reynolds and the Mach numbers are sufficiently small, then time periodic parallel flows are asymptotically stable under perturbations small in some Sobolev space on the layer $\tilde \Omega$. Furthermore, it was shown that the asymptotic leading part of the perturbation is given by a product of a time periodic function and a solution of an $n-1$ dimensional linear  
heat equation in the case $n = 3$, and by a product of a time periodic function and a solution of a one-dimensional viscous Burgers equation in the case $n=2$; the hyperbolic aspect of the perturbation decays faster. 
(See \cite{Brezina3} and references therein for the mathematical analysis of the stability of parallel flows in compressible fluids.)

On the other hand, in reality, the external force $\tilde G_{para} = \trans(\tilde g^{1}(\tilde x_{n}, \tilde t), 0, \cdots, 0)$ often undergoes a perturbation in $\tilde x'$ variable. Under such a situation, the external force depends not only $(\tilde x_{n},\tilde t)$ but also $\tilde x'$, so the time periodic parallel flow is no longer a solution of (1.1)-(1.2) since $\tilde G$ depends on $\tilde x'$. 
In this paper, we thus consider the situation where the external force $\tilde G$ periodically depends on $\tilde x'$ variable as described in \eqref{force}. 
Under such a situation, as was mentioned above, if $\tilde{G}$ is sufficiently small, the system (\ref{den1})-(\ref{vel1}) with (\ref{cond}) has a space-time periodic state $\tilde{u}_p=\trans(\tilde{\rho}_p, \tilde{v}_p)$. 
We shall establish the results on the spectral properties of the linearized evolution operator around $\tilde u_{p}$ which suggest that the asymptotic leading part of the perturbation of $\tilde u_{p}$ exhibits diffusive behaviors similar to those in the case of parallel flows in \cite{Brezina3} if the Reynolds and the Mach numbers are sufficiently small.   

We briefly explain our main results of this paper. 
After introducing suitable non-dimensional variables, the equations for the perturbation 
$u(t)= \trans(\phi,w)= \trans(\gamma^2 (\rho-\rho_p),v-v_p)$ takes the following form: 
\begin{gather}\label{den2}
\begin{split} 
\partial_t \phi + \mathrm{div}(\phi v_p)+\gamma^2 \mathrm{div}(\rho_p w)=f(u), 
\end{split} \\ \label{vel2}
\begin{split}  
\partial_t w &-\frac{\nu }{\rho_p}\Delta w-\frac{\tilde{\nu} }{\rho_p}\nabla \mathrm{div}v + \nabla \big(\frac{p'(\rho_p)}{\gamma^2 \rho_p}\phi \big) \\ 
&+ \frac{1}{\gamma^2 \rho_p^2}(\nu \Delta v_p  
+ \tilde{\nu}\nabla \mathrm{div}v_p)\phi + v_p\cdot \nabla w +w\cdot \nabla v_p=g(u), 
\end{split} 
\end{gather}
on $\Omega= \mathbb{R}^{n-1} \times (0,1)$, and 
\begin{gather}
w|_{\partial \Omega}=0, \label{boundary} \\
(\phi,w)|_{t=0}=(\phi_0,w_0). \label{initial}
\end{gather}
Here $u_p=\trans(\phi_p, w_p)$ denotes the non-dimensionalization of $\tilde{u}_p=\trans(\tilde{\phi}_p, \tilde{w}_p)$, and $\nu$, $\tilde{\nu}$ and $\gamma$ are non-dimensional parameters. 
The terms $f(u)$ and $g(u)$ are non-linear terms given by 
\begin{align*}
f(u)=&\mathrm{div}(\phi w),\\
g(u)=&-w\cdot \nabla w -\frac{\phi}{\rho_p(\gamma^2 \rho_p + \phi)}\big(\nu\Delta w+\tilde{\nu}\nabla \mathrm{div}w -\frac{\nu \phi}{\gamma^2 \rho_p} \Delta v_p -\frac{\tilde{\nu} \phi}{\gamma^2 \rho_p}\nabla \mathrm{div}v_p   \big)\\
&+\frac{\phi}{\gamma^2 \rho_p^2}\nabla\big(p^{(1)}\big(\rho_p,\frac{\phi}{\gamma^2}\big)\frac{\phi}{\gamma^2}\big)\\
&+\frac{\phi^2}{\gamma^2 \rho_p^2(\gamma^2 \rho_p + \phi)}\nabla\big(p\big(\rho_p+\frac{\phi}{\gamma^2}\big)\big)+\frac{1}{\rho_p}\nabla\big(p^{(2)}\big(\rho_p,\frac{\phi}{\gamma^2}\big)\frac{\phi^2}{\gamma^4}\big),
\end{align*}
where 
\begin{align*}
&p^{(1)}(\rho_p,\phi)=\int_{0}^{1}p'(\rho_p + \theta \phi)d\theta ,\\
&p^{(2)}(\rho_p,\phi)=\int_{0}^{1}(1-\theta)p''(\rho_p + \theta \phi)d\theta .
\end{align*}

We consider the linearized problem for (\ref{den2})-(\ref{initial}) which can be written as 
\begin{equation} \label{linearized}
\partial_t u+ L(t)u = 0, \  u|_{t=s}=u_0 
\end{equation}
on $\Omega= \mathbb{R}^{n-1} \times (0,1)$, where $u(t)=\trans(\phi(t), w(t))\in D(L(t))$; and $L(t)$ is the operator on $L^2(\Omega)\times L^2(\Omega)$ of the form  
\begin{equation*} 
\begin{split}
L(t)=&
\begin{pmatrix}
\mathrm{div}( v_p(t) \cdot) & \gamma^2 \mathrm{div}(\rho_p(t)\cdot)\\
\nabla \big(\frac{p'(\rho_p(t))}{\gamma^2 \rho_p(t)}\cdot\big) & -\frac{\nu}{\rho_p(t)}\Delta-\frac{\tilde{\nu}}{\rho_p(t)}\nabla \mathrm{div}
\end{pmatrix} \\
&+\begin{pmatrix} 
0 & 0\\
\frac{1}{\gamma^2 \rho_p^2(t)}(\nu\Delta v_p(t)+\tilde{\nu}\nabla \mathrm{div} \  v_p(t) ) & v_p(t) \cdot \nabla+\trans(\nabla v_p(t) )
\end{pmatrix}
\end{split}
\end{equation*}
with domain 
\begin{equation*}
D(L(t))=\{u=\trans(\phi,w)\in L^2(\Omega) ;w \in H_0^1(\Omega),L(t)u \in L^2(\Omega)\}.
\end{equation*}

We denote by $\mathbf{U}(t,s)$ the solution operator for (\ref{linearized}). 
Since $L(t)$ has spatially periodic coefficients, the Bloch transform is useful to study the spectral properties of $\mathbf{U}(t,s)$. If we apply the Bloch transform to $\mathbf{U}(t,s)$, we have a family $\{\mathbf{U}_{\eta'}(t,s)\}_{\eta'\in Q^*}$ of solution operators, where $Q^*=\Pi_{i=1}^{n-1}[-\frac{\alpha_i}{2}, \frac{\alpha_i}{2})$; and each $\mathbf{U}_{\eta'}(t,s)$ is the solution operator for the problem 
\begin{equation}
\partial_t u+ L_{\eta'}(t)u = 0, \  u|_{t=s}=u_0 \label{linearized eta}
\end{equation}
on 
\begin{equation*}
\Omega_{per}= \Pi^{n-1}_{j=1}\mathbb{T}_{\frac{2\pi}{\alpha_j}} \times (0,1). 
\end{equation*}
Here $\mathbb{T}_a=\mathbb{R}/a\mathbb{Z}$; and $L_{\eta'}(t)$ is an operator acting on functions on $\Omega_{per}$ which takes the form 
\begin{align*}
L_{\eta'}(t)=&
\begin{pmatrix}
\nabla_{\eta'}\cdot( v_p(t) \cdot) & \gamma^2 \nabla_{\eta'}\cdot(\rho_p(t)\cdot)\\
\nabla_{\eta'} \big(\frac{p'(\rho_p(t))}{\gamma^2 \rho_p(t)}\cdot\big) & -\frac{\nu}{\rho_p(t)}\Delta_{\eta'}-\frac{\tilde{\nu}}{\rho_p(t)}\nabla_{\eta'} \trans\nabla_{\eta'}
\end{pmatrix} \\
&+\begin{pmatrix} 
0 & 0\\
\frac{1}{\gamma^2 \rho_p^2(t)}(\nu\Delta v_p(t)+\tilde{\nu}\nabla \mathrm{div} \  v_p(t)) & v_p(t) \cdot \nabla_{\eta'}+\trans(\nabla v_p(t) )
\end{pmatrix},
\end{align*} 
where $\nabla_{\eta'}$ and $\Delta_{\eta'}$ are defined by 
\begin{equation*}
\nabla_{\eta'} = \nabla +i \tilde{\eta}', \hspace{8mm} \Delta_{\eta'} = \nabla_{\eta'}\cdot \nabla_{\eta'}, 
\end{equation*}
with $\tilde{\eta}'= \trans(\eta',0)\in \mathbb{R}^n$.

As in \cite{Brezina2, Brezina1}, we investigate the spectral properties of $\mathbf{U}_{\eta'}(t,s)$ by the Floquet theory. 
Since $L_{\eta'}(t+1)=L_{\eta'}(t)$ for all $t$, the large time behavior of $\mathbf{U}_{\eta'}(t,s)$ is controlled by the spectrum of the monodromy operator $\mathbf{U}_{\eta'}(1,0)$. 
We thus consider the spectral properties of $B_{\eta'}$, 
where $B_{\eta'}$ is an operator on the time periodic function space $X=L^2(\mathbb{T}_1; L^2(\Omega_{per})\times L^2(\Omega_{per}))$ defined by 
\begin{align*}
D(B_{\eta'})&=\{ u=\trans(\phi,w)\in X; 
w \in H^1(\mathbb{T}_1;H^{-1}(\Omega_{per}) ) \cap L^2(\mathbb{T}_1; H^1_0(\Omega_{per})), \\
&\hspace{8mm}\phi \in C(\mathbb{T}_1;L^2(\Omega_{per})), 
B_{\eta'}u \in X \}, \\
B_{\eta'}u&=\partial_tu+L_{\eta'}u, \hspace{4mm} u=\trans(\phi, w)\in D(B_{\eta'}).
\end{align*} 
The spectrum of $-B_{\eta'}$ gives Floquet exponents of the problem (\ref{linearized eta}). 

Our main results are summarized as follows. 
If the external force $G$ and Bloch parameter $\eta'$ are sufficiently small, then 
\begin{equation}\label{sec1main1}
\sigma(-B_{\eta'}) \cap\left\{\lambda\in\mathbb{C}; \mathrm{Re}\lambda> -\frac{\beta_0}{4},\, |\lambda-2\pi i k|\leq \frac{\pi}{4},\, k\in \mathbb{Z}\right\} =\{\lambda_{\eta',k}: k\in\mathbb{Z}\}. 
\end{equation}
Here $\beta_{0}$ is a positive constant; and $\lambda_{\eta',k}$ is a simple eigenvalue of $-B_{\eta'}$ satisfying 
\begin{equation}\label{sec1main2}
\lambda_{\eta',k}=2\pi ik- i\sum_{j=1}^{n-1}a_j\eta_j-\sum_{j,k=1}^{n-1}a_{jk}\eta_j \eta_k  + O(|\eta'|^3) \hspace{4mm}(\eta'\rightarrow 0)
\end{equation}
with some constants $a_j, a_{jk} \in \mathbb{R}$, where $(a_{jk})_{1\leq j,k\leq n-1}$ is positive definite. 
It follows from \eqref{sec1main1} that the spectrum of $\mathbf{U}_{\eta'}(1,0)$ with $|\eta'| \ll 1$ satisfies  
\[
\sigma(\mathbf{U}_{\eta'}(1,0) ) \subset \{\lambda; |\lambda| \le e^{-\frac{\beta_{0}}{4}}\} \cup \{e^{\lambda_{\eta',0}} \}, 
\]
where $e^{\lambda_{\eta',0}}$ is a simple eigenvalue of $\mathbf{U}_{\eta'}(1,0)$. This yields the asymptotic behavior $\mathbf{U}_{\eta'}(m,0) \, \sim \, e^{m \lambda_{\eta',0}}$ as $m \to \infty$, which, together with  \eqref{sec1main2}, would imply that $\mathbf{U}(t,0)$ would behave diffusively as $t \to \infty$. 
We also establish the boundedness of the eigenprojection for the eigenvalue $\lambda_{\eta', 0}$ with $|\eta'| \ll 1$ which is needed in the analysis of the nonlinear problem. 

Our results, in fact, will yield the following diffusive behavior of a part of the solution operator $\mathbf{U}(t,s)$ as $t - s \to \infty$ that includes the space-time periodic nature of the problem. In a similar manner to \cite{Brezina3,Brezina2, Brezina1}, based on \eqref{sec1main1}, \eqref{sec1main2} and the Floquet theory, one can show that there exist a bounded projection $\mathbb{P}(t)$ on $L^2(\Omega)\times L^2(\Omega)$ such that $\mathbb{P}(t+1)= \mathbb{P}(t)$ and the following estimates hold:
\begin{gather*}
\| \mathbb{P}(t)\mathbf{U}(t,s)u_0\|_{L^2(\Omega)\times L^2(\Omega)}\leq C(1+t-s)^{-\frac{n-1}{4}}\| u_0\|_{L^1(\Omega)\times L^1(\Omega)}, \\ 
\begin{split}
&\| \mathbb{P}(t)\mathbf{U}(t,s)u_0 -u^{(0)}(t)\mathcal{H}(t-s)\sigma_0 \|_{L^2(\Omega)\times L^2(\Omega)} \\
&\leq C(t-s)^{-\frac{n-1}{2}(\frac{1}{p}-\frac{1}{2})-\frac{1}{2}}\| u_0\|_{L^p(\Omega)\times L^p(\Omega)} \hspace{4mm}(1\leq p\leq 2). 
\end{split}
\end{gather*}
Here $u^{(0)}=u^{(0)}(x',x_n,t)$ is some function $\frac{2\pi}{\alpha_i}$-periodic in $x_i~(i=1, \cdots, n-1)$ and $1$-periodic in $t$ and $\mathcal{H}(t)\sigma_0$ is a solution of the linear heat equation
\begin{equation*}
\begin{cases}
\partial_t \sigma +\sum_{j=1}^{n-1}a_j\partial_{x_j} \sigma - \sum_{j,k=1}^{n-1}a_{jk}\partial_{x_j}\partial_{x_k}\sigma=0,\\
\sigma|_{t=0}= \sigma_0.
\end{cases}
\end{equation*}

The main difference in the analysis of this paper to the case of the parallel flow in \cite{Brezina2, Brezina1} is as follows. In the case of the parallel flow, by the Fourier transform in $x' \in \mathbb{R}^{n-1}$, the spectral analysis for the linearized problem is reduced to the one for a {\it one-dimensional problem} on the interval $(0,1)$ with a parameter of the Fourier variable $\xi' \in \mathbb{R}^{n-1}$; and the one-dimensional aspect of the reduced problem was essentially used in the analysis in \cite{Brezina2, Brezina1}, e.g., to obtain a regularity estimate of time-periodic eigenfunctions for the Floquet exponents of the operator corresponding to $-B_{\eta'}$ with $|\eta'|\ll1$. On the other hand, the Bloch transformed problem \eqref{linearized eta} is a {\it multi-dimensional problem}, i.e., a problem on $\Omega_{per}$, which requires approaches different to those in \cite{Brezina2, Brezina1}, e.g., we specify the eigenspace for the eigenvalue $0$ of $-B_0$ and construct time periodic eigenfunctions of $-B_{\eta'}$ with $|\eta'|\ll1$ in a higher order Sobolev space, based on the energy methods in \cite{Iooss-Padula,M-N} and the argument to construct time periodic solutions in \cite{Valli}. 

We also mention that dissipative systems on infinite layers and cylindrical domains often provide space-time periodic patterns (cf., \cite{Doelman-Sandstede-Scheel-Schneider, Schneider}). The analysis of this paper is thus a preparatory study of the dynamics around space-time periodic patterns of the viscous compressible system \eqref{den1}-\eqref{vel1}. 

This paper is organized as follows. In Section 2, we transform the equations (\ref{den1})-(\ref{vel1}) into a non-dimensional form and introduce basic notation that is used throughout the paper. In Section 3, we first state the existence of a space-time periodic state and then state the main results of this paper. Section 4 is devoted to the proof of the main results. 
In Section 5, we give a proof of the existence of a space-time periodic state.

\section{Preliminaries}
In this section, we transform (\ref{den1})-(\ref{vel1}) into a non-dimensional form and introduce some function spaces and notations which are used throughout the paper. 

\vspace{1ex}

We rewrite the problem into the non-dimensional form.
We introduce the following non-dimensional variables:
\begin{gather*}
\tilde{x}=dx, \hspace{4mm}\tilde{t}=Tt, \hspace{4mm} \tilde{\rho}=\rho_*\rho, \hspace{4mm} \tilde{v}=\frac{d}{T}v, \hspace{4mm} \tilde{p}=\rho_* \tilde{p}'(\rho_*)p, \hspace{4mm}\tilde{G}=[\tilde{G}]_{3,T,\tilde{\Omega}_{per}}G,
\end{gather*}
where 
\begin{align*}
[\tilde{G}]_{3,T,\tilde{\Omega}_{per}}&=\left(\sum_{j=0}^2T^{2j-1}d^{2(3-2j)-n}\int_0^T\|\partial_{\tilde{t}}^j\tilde{G}\|^2_{H^{3-2j}(\tilde{\Omega}_{per})}\,d\tilde{t}\right)^{\frac{1}{2}}. 
\end{align*}
Here $\tilde{\Omega}_{per}=\Pi^{n-1}_{j=1}\mathbb{T}_{\frac{2\pi}{\tilde{\alpha}_j}} \times (0,d)$; and $\|\cdot\|_{H^m(\tilde{\Omega}_{per})}$ denotes the usual $H^m$-norm over $\tilde{\Omega}_{per}$ (whose definition is given below).

Under this change of variables, the domain $\tilde{\Omega}$ is transformed into 
$$\Omega = \mathbb{R}^n \times (0,1).$$ 
The equations (\ref{den1})-(\ref{vel1}) are rewritten as 
\begin{gather} \label{dim1}
\partial_t \rho + \mathrm{div}(\rho v)=0, \\
\rho(\partial_t v+v \cdot \nabla v)-\nu \Delta v-\tilde{\nu}\nabla \mathrm{div}v +\gamma^2 \nabla p(\rho)=S\rho G. \label{dim2}
\end{gather} 
Here $\nu$, $\tilde{\nu}$, $\gamma$ and $S$ are non-dimensional parameters defined by 
\begin{equation*}
\nu=\frac{\mu T}{\rho_*d^2},\hspace{4mm} \tilde{\nu}=(\mu+\mu')\frac{T}{\rho_*d^2},\hspace{4mm} 
\gamma= \frac{T}{d}\sqrt{\tilde{p}'(\rho_*)},\hspace{4mm} S=\frac{T^2}{d}[\tilde{G}]_{3,T,\tilde{\Omega}_{per}}.
\end{equation*}
We note that 
\begin{equation*}
p'(1)=1 \quad \text{and}\quad [G]_{3,1,\Omega_{per}} =1.
\end{equation*}
Furthermore, due to the assumption (\ref{assump mu}), we have 
\begin{equation*}
\nu\leq \nu+\tilde{\nu}\leq (2+\mu_*)\nu. 
\end{equation*}

The boundary and initial conditions (\ref{cond})-(\ref{cond1}) are transformed into 
\begin{gather} \label{dim3}
v|_{x_n =0,1}=0 , \\ \label{dim4}
(\rho,v)|_{t=0}= (\rho_0,v_0).
\end{gather}

We next introduce notation used throughout this paper. 
Let $D$ be a domain. We denote by $L^p(D)$ $(1\leq p \leq \infty)$ the usual Lebesgue space on $D$ and its norm is denoted by $\|\cdot\|_{L^p(D)}$. Let $m$ be a nonnegative integer. $H^m(D)$ denotes the $m$-th order $L^2$-Sobolev space on $D$ and its norm denoted by $\|\cdot\|_{H^m(D)}$. $C_0^m(D)$ is defined as the set of $C^m$-functions having compact supports in $D$. Furthermore,  we denote by $H_0^m(D)$ the completion of $C_0^{\infty}(D)$ in $H^m(D)$ and the dual space of $H_0^m(D)$ is denoted by $H^{-m}(D)$.

We simply write the set of all vector fields $w=\trans(w_1,\cdots,w_n)$ on $D$ as $w_j \in L^p(D)$ (resp., $H^m(D)$) and its norm is denoted by $\|\cdot\|_{L^p(D)}$ (resp., $\|\cdot\|_{H^m(D)}$). For $u=\trans(\phi,w)$ with $\phi \in H^k(D)$ and $w \in H^m(D)$, we define $\|u\|_{H^k(D)\times H^m(D)}=(\|\phi\|_{H^k(D)}^2+ \|w\|_ {H^m(D)}^2)^{\frac{1}{2}}$. 
When $k=m$, we simply write $\|u\|_{H^k(D)\times H^k(D)}=\|u\|_{H^k(D)}$.

We set 
\begin{equation*}
\nom{f(t)}{}{k}= \big(\sum_{j=0}^{\left[\frac{k}{2}\right]}\|\partial_t^j f(t)\|_{H^{k-2j}(\Omega_{per})}^2\big)^\frac{1}{2}, 
\end{equation*}
where $[k]$ is the largest integer smaller than or equal to $k$.

The inner product of $L^2$ is defined as 
\begin{equation*}
(f,g)=\int_{\Omega_{per}}f(x) \overline{g(x)}\, dx
\end{equation*}
for $f,g \in L^2(\Omega_{per})$. Here, $\overline{g}$ denotes the complex conjugate of $g$. 
Moreover, the mean value of $f=f(x)$ and $g=g(x,t)$ over $\Omega_{per}$ and $\Omega_{per}\times\mathbb{T}_1$ is written as 
\begin{equation*}
\langle f \rangle=\int_{\Omega_{per}}f(x)\, dx \quad \text{and}\quad  \mean{g}=\int^1_0\langle g(t)\rangle\, dt,
\end{equation*}
respectively. 
We next introduce a weighted inner product:
\begin{equation*}
\inn{u_1}{u_2}=\int_0^1 \langle u_1(t),u_2(t)\rangle_t\, dt
\end{equation*}
for $u_j=\trans(\phi_j,w_j)$ $j=1,2$, where 
\begin{equation*}
\langle u_1(t), u_2(t)\rangle_t=\int_{\Omega_{per}}\phi_1(t) \overline{\phi_2(t)} \frac{p'(\rho_p(t))}{\gamma^4 \rho_p(t)} + w_1(t) \cdot\overline{w_2(t)} \rho_p(t)\, dx. 
\end{equation*}
Here $\rho_p$ denotes the density of the space-time periodic state given in Proposition \ref{timeper} below. 
By Proposition \ref{timeper}, we see that $\rho_p\geq\underline{\rho}$ on $\Omega_{per}\times\mathbb{T}_1$ for a positive constant $\underline{\rho}$ and that $|\rho_p(x,t)-1|\leq \frac{1}{2}$ and $|p'(\rho_p(x,t))-1|\leq\frac{1}{2}$ for all $(x,t)\in\Omega_{per}\times \mathbb{T}_1$. 
Therefore, $\inn{u_1}{u_2}$ defines an inner product. 

We finally define $L_{*}^2(\Omega_{per})$ and $H_{*}^m(\Omega_{per})$ by 
\begin{gather*}
L_*^2(\Omega_{per}):=\{\phi \in  L^2(\Omega_{per});\langle \phi \rangle = 0 \}
\end{gather*}
and 
\begin{gather*}
H_*^m(\Omega_{per})= H^m(\Omega_{per}) \cap L_*^2(\Omega_{per}),
\end{gather*} 
respectively.

We next introduce the Bogovskii lemma \cite{Bogovskii,G}.
\begin{lem}[\cite{Bogovskii,G}]\label{lem Bogovskii}
There exist a bounded operator $\mathcal{B}:L_*^2(\Omega_{per}) \rightarrow H_0^1(\Omega_{per})$ such that for any $f \in L_*^2(\Omega_{per})$, 
\begin{gather*}
\mathrm{div} \mathcal{B}f = f, \\
\|\nabla \mathcal{B}f\|_{L^2(\Omega_{per})}  \leq C \|f\|_{L^2(\Omega_{per})}, 
\end{gather*}
where $C$ is a positive constant depending only on $\Omega_{per}$. 
Furthermore, if $f=\mathrm{div}g$ with $g=\trans(g^1,\cdots, g^n)$ satisfying $g\in H^1_0(\Omega_{per})$, then 
\begin{gather*}
\| \mathcal{B}(\mathrm{div}g)\|_{L^2(\Omega_{per})}  \leq C \|g\|_{L^2(\Omega_{per})}.
\end{gather*}
\end{lem}

In terms of the Bogovskii operator $\mathcal{B}$, we introduce the following inner product on $L^2_*(\Omega_{per})\times L^2(\Omega_{per})$. 
For each $t\in\mathbb{T}_1$, we define $(\!(u_1,u_2)\!)_t$ by
\begin{equation*}
(\!(u_1,u_2)\!)_t =\langle u_1,u_2\rangle_t - \delta [(w_1,\mathcal{B}\phi_2) + (\mathcal{B}\phi_1,w_2)], 
\end{equation*} 
where $\delta$ is a positive constant. 
One can see that there exists a positive constant $C$ such that if $0<\delta <\frac{1}{2C\gamma}$, then $(\!(\cdot, \cdot)\!)_t$ defines an inner product satisfying 
\begin{equation*}
\frac{1}{2}\|u\|^2_{L^2(\Omega_{per}),\gamma}\leq (\!(u,u)\!)_t\leq \frac{3}{2}\|u\|^2_{L^2(\Omega_{per}),\gamma}, 
\end{equation*} 
where 
\begin{equation*}
\|u\|^2_{L^2(\Omega_{per}),\gamma}=\frac{1}{\gamma^2}\|\phi\|^2_{L^2(\Omega_{per})}+\|w\|^2_{L^2(\Omega_{per})}.
\end{equation*}

We next introduce the Bloch transform. 
Let $\mathscr{S}\,(\mathbb{R}^{n-1})$ be the Schwartz space on $\mathbb{R}^{n-1}$. 
We define the Bloch transform $\mathcal{T}$ by  
\begin{align*}
(\mathcal{T}\varphi)(x',\eta')&=\frac{1}{(2\pi)^{\frac{n-1}{2}}|Q|^{\frac{1}{2}}} \sum_{k_1,...,k_{n-1}\in \mathbb{Z}^{n-1}} \hat{\varphi}(\eta' + \sum_{j=1}^{n-1}k_j\alpha_j \mathbf{e}_j')\textstyle{e^{i\sum_{j=1}^{n-1}k_j\alpha_j x'}}\\
&=\frac{1}{|Q^*|^{\frac{1}{2}}}\sum_{l_1,...l_{n-1}\in \mathbb{Z}^{n-1}} \varphi(x' + \sum_{j=1}^{n-1}l_j\frac{2\pi}{\alpha_j} \mathbf{e}_j')\textstyle{e^{-i\eta'\cdot (x'+\sum_{j=1}^{n-1}l_j\frac{2\pi}{\alpha_j} \mathbf{e}_j')}}
\end{align*}
for $\varphi\in\mathscr{S}(\mathbb{R}^{n-1})$, 
where $\hat{\varphi}$ denotes the Fourier transform of $\varphi$: 
$$\hat{\varphi}(\xi')=\int_{\mathbb{R}^{n-1}} \varphi(x')e^{-i \xi'\cdot x'}\,dx';$$
and 
$$Q=\prod_{i=1}^{n-1}\big[-\frac{\pi}{\alpha_i},\frac{\pi}{\alpha_i}\big), \quad  Q^*=\prod_{i=1}^{n-1}\big[-\frac{\alpha_i}{2},\frac{\alpha_i}{2}\big).$$

Let $\varphi(x',\eta')$ be in $C^{\infty}(\mathbb{R}^{n-1}\times \mathbb{R}^{n-1})$ such that $\varphi(x',\eta')$ is $Q$-periodic in $x'$ and $\varphi(x',\eta')e^{i\eta'\cdot
 x'}$ is $Q^*$-periodic in $\eta'$. We define $(\mathcal{S}\varphi)(x')$ by 
\begin{equation*}
(\mathcal{S}\varphi)(x')=\frac{1}{|Q^*|^\frac{1}{2}}\int_{Q^*}\varphi(x',\eta')e^{i\eta'\cdot x'}\,d\eta',
\end{equation*}
where $x' \in \mathbb{R}^{n-1}$. 
Note that $\varphi(x',\eta'+\alpha_j\mathbf{e}_j')=\varphi(x',\eta')e^{-i\alpha_j \mathbf{e}_j'\cdot x'}$. 

The operators $\mathcal{T}$ and $\mathcal{S}$ have the following properties. See, e.g., \cite{R-S,Schneider} for the details. 
\begin{prop} \label{Bloch}
$(\mathrm{i})$ $(\mathcal{T}\varphi)(x',\eta')$ is $Q$-periodic in $x'$ and $(\mathcal{T}\varphi)(x',\eta')e^{i\eta'.x'}$ is $Q^*$-periodic in $\eta'$. 

\vspace{1ex}
$(\mathrm{ii})$ $\mathcal{T}$ is uniquely extended to an isometric operator from $L^2(\mathbb{R}^{n-1})$ to $L^2(Q^*;L^2(Q))$. 

\vspace{1ex}
$(\mathrm{iii})$ $\mathcal{S}$ is the inverse operator of $\mathcal{T}$. 

\vspace{1ex}
$(\mathrm{iv})$ Let $\varphi$ be $Q$-periodic in $x'$. Then it holds that $\mathcal{T}(\psi \varphi)=\psi \mathcal{T}\varphi$. 

\vspace{1ex}
$(\mathrm{v})$ $\mathcal{T}(\partial_{x_j}\varphi)=(\partial_{x_j}+i\eta_j)\mathcal{T}\varphi$ and $\mathcal{T}$ defines  an isomorphism from $H^m(\mathbb{R}^{n-1})$ to $L^2(Q^*;H^m(Q))$. 
\end{prop}

\section{Main Results}
In this section, we state the main results of this paper. 
We first state the existence of the space-time periodic state of (\ref{dim1})-(\ref{dim3}). 
We consider the time periodic problem for 
\begin{gather} \label{prob time1}
\partial_t \rho + \mathrm{div}(\rho v)=0, \\
\rho(\partial_t v+v \cdot \nabla v)-\nu \Delta v-\tilde{\nu}\nabla \mathrm{div}v +\gamma^2 \nabla p(\rho)=S\rho G \label{prob time2}
\end{gather} 
in $\Omega_{per}$ under the boundary condition 
\begin{gather} 
v|_{x_n =0,1}=0. \label{boundary time}
\end{gather}

\begin{prop}\label{timeper}
Let $G\in \cap_{j=0}^2C^j(\mathbb{T}_1; H^{3-2j}(\Omega_{per}))$ with $[G]_{3,1,\Omega_{per}}=1$.
There exist positive constants $\nu_0$, $\gamma_0$, $\varepsilon_{0}$ and $a$ such that if $\frac{\nu^2}{\nu+\tilde{\nu}} \geq \nu_0$, $\frac{\gamma^2}{\nu + \tilde{\nu}} \geq \gamma_0$ and $S\leq \varepsilon_{0}\frac{\nu^2}{\gamma^2\sqrt{\nu+\tilde{\nu}}}\sqrt{1-e^{-a\frac{\nu+\tilde{\nu}}{\gamma^2}}}$, then the following assertions hold. 
There exist a space-time periodic solution $u_p=\trans(\rho_p,v_p)=\trans(1+\phi_p,v_p)\in \cap_{j=0}^2 C^j(\mathbb{T}_1; H^{4-2j}(\Omega_{per})\times H^{4-2j}(\Omega_{per}))\cap H^j(\mathbb{T}_1;H^{4-2j}(\Omega_{per})\times H^{5-2j}(\Omega_{per}))$ of problem $(\ref{prob time1})$-$(\ref{boundary time})$ satisfying $\langle\rho_p(t)\rangle=1$ for each $t\in\mathbb{T}_1$  and $\rho_p=\rho_p(x,t)\geq \underline{\rho}$ for a positive constant $\underline{\rho}$. 
Furthermore, $u_p$ satisfies the following estimates 
\begin{gather} \label{timeperiodic}
\gamma^2\nom{\phi_p (t)}{2}{4} + \nom{v_p (t)}{2}{4} \leq C\frac{\nu^2}{\gamma^4}, \\
\int_{0}^{1}\frac{\gamma^4}{\nu+\tilde{\nu}}\nom{\nabla \phi_p (s)}{2}{3} +(\nu+\tilde{\nu})\|\partial_t^2\phi\|_{L^2(\Omega_{per})}^2 + \frac{\nu^2}{\nu+\tilde{\nu}} \nom{v_p (s)}{2}{5}\,dt
\leq C\frac{\nu^2}{\gamma^4}, \label{timeperiodic2}
\end{gather}
where $C$ is a positive constant independent of $\nu$, $\tilde{\nu}$, $\gamma$ and $S$.
\end{prop}

\begin{rem}
$(\mathrm{i})$ Since $\|\phi_p\|_{L^\infty(\Omega_{per})}\ll1$ if $\frac{\gamma^2}{\nu+\tilde{\nu}}\gg1$ and $\frac{\nu^2}{\nu+\tilde{\nu}}\gg1$, we have $\rho_p\sim 1$, and therefore $p'(\rho_p)\sim1$. \\
$(\mathrm{ii})$ If $\frac{\gamma^2}{\nu+\tilde{\nu}}\gg1$, then the assumption on $S$ in Proposition $\ref{timeper}$ implies 
\begin{equation*}
S\leq \varepsilon_0\sqrt{a}\frac{\nu^2}{\gamma^3}.
\end{equation*} 
\end{rem}

The proof of Proposition \ref{timeper} is essentially the same as that given in \cite{Valli}. 
Since we solve the time periodic problem in $H^4(\Omega_{per})\times H^5(\Omega_{per})$ and we need to know the dependence of the estimates on the parameters $\nu$, $\tilde{\nu}$ and $\gamma$, we will give an outline of the proof of Proposition \ref{timeper} in Section \ref{sec5} below.

Our main result is concerned with the spectrum of the linearized solution operator $\mathbf{U}(t,s)$ around the space-time periodic solution $u_p=\trans(\phi_p, v_p)$. 

As was mentioned in the introduction, we apply Bloch transform to (\ref{linearized}). 
By Proposition \ref{Bloch}, we then obtain (\ref{linearized eta}) and consider the spectrum of $-B_{\eta'}$ to obtain the Floquet exponents of (\ref{linearized eta}) for $|\eta'|\ll1$.
\begin{thm} \label{critical eigen}
There exist positive constants $\nu_0$, $\gamma_0$, $\varepsilon_0$ and $a$  such that if $\frac{\nu^2}{\nu+\tilde{\nu}} \geq \nu_0$, $\frac{\gamma^2}{\nu+ \tilde{\nu}} \geq \gamma_0$ and $S\leq \varepsilon_{0}\frac{\nu^2}{\gamma^2\sqrt{\nu+\tilde{\nu}}}\sqrt{1-e^{-a\frac{\nu+\tilde{\nu}}{\gamma^2}}}$, then the following assertions hold. 

{\rm (i)} There exists a positive constant $r_0=r_0(\nu, \tilde{\nu}, \gamma, \nu_0, \gamma_0)$ such that if $|\eta'|\leq r_0$, then  
\begin{gather*}
\Sigma:=\left\{\lambda\in\mathbb{C}; \mathrm{Re}\lambda> -\frac{\beta_0}{4},\, |\lambda-2\pi i k|\geq \frac{\pi}{4},\, k\in \mathbb{Z}\right\}\subset \rho(-B_{\eta'}), 
\end{gather*}
and for $\lambda\in\Sigma$ 
\begin{gather*}
\|(\lambda +B_{\eta'})^{-1}F\|_{L^2(\mathbb{T}_1;L^2(\Omega_{per})\times H^1(\Omega_{per}))}\leq C\| F\|_X.
\end{gather*}

{\rm (ii)} If $|\eta'|\leq r_0$, then 
\begin{gather*}
\sigma(-B_{\eta'}) \cap\left\{\lambda\in\mathbb{C}; \mathrm{Re}\lambda> -\frac{\beta_0}{4},\, |\lambda-2\pi i k|\leq \frac{\pi}{4},\, k\in \mathbb{Z}\right\} =\{\lambda_{\eta',k}: k\in\mathbb{Z}\},
\end{gather*} 
where $\lambda_{\eta',k}$ is a simple eigenvalue that satisfies 
\begin{equation*}
\lambda_{\eta',k}=2\pi ik-i\sum_{j=1}^{n-1}a_j\eta_j - \sum_{j,k=1}^{n-1}a_{jk}\eta_j \eta_k+ O(|\eta'|^3) \hspace{4mm} (\eta' \rightarrow 0)
\end{equation*}
with some constants $a_j,a_{jk} \in \mathbb{R}$ satisfying 
\begin{equation*}
\sum_{j,k=1}^{n-1}a_{jk}\xi_j \xi_k \geq \frac{\kappa_0 \gamma ^2}{\nu}|\xi'|^2,
\end{equation*}
for  all $\xi'=\trans(\xi_1,\cdots,\xi_{n-1}) \in \mathbb{R}^{n-1}$ 
and some positive constant $\kappa_0$.
As a consequence 
\begin{equation}
\mathrm{Re}\lambda_{\eta',k}\leq -\frac{\kappa_0 \gamma ^2}{2 \nu}|\eta'|^2 .\label{eigen}
\end{equation}
\end{thm}

We next consider eigenfunctions for eigenvalues $\lambda_{\eta',0}$. 
We introduce the adjoint operator $B_{\eta'}^*$ defined by 
\begin{gather*}
\begin{split} 
 D(B_{\eta'}^*)=\{& u=\trans(\phi,w)\in X; 
w \in H^1(\mathbb{T}_1;H^{-1}(\Omega_{per}))\cap L^2(\mathbb{T}_1;H^1_0(\Omega_{per})), \\
&\phi \in C(\mathbb{T}_1;L^2(\Omega_{per})), 
B^*_{\eta'}u \in X \}, 
\end{split} \\
B_{\eta'}^*u=-\partial_t^* u + L_{\eta'}^*(\cdot)u, \hspace{4mm}u=\trans(\phi,w)\in D(B_{\eta'}^*), 
\end{gather*}
where 
\begin{gather*}
\partial_t^* = 
\begin{pmatrix}
\frac{\gamma^4 \rho_p}{p'(\rho_p)} \partial_t(\frac{p'(\rho_p)}{\gamma^4 \rho_p} \cdot) \\ \frac{1}{\rho_p}\partial_t(\rho_p \cdot)
\end{pmatrix}, \\
\begin{split}
L_{\eta'}^*=&
\begin{pmatrix}
-v_p\cdot \nabla_{\eta'}\big(\frac{p'(\rho_p)}{ \rho_p}\cdot \big)\frac{\rho_p}{p'(\rho_p)}& -\gamma^2\trans \nabla_{\eta'}\cdot (\rho_p \cdot)\\
-\nabla_{\eta'}\big(\frac{p'(\rho_p)}{\gamma^2 \rho_p}\cdot\big) & -\frac{\nu}{\rho_p}\Delta_{\eta'}-\frac{\tilde{\nu}}{\rho_p}\nabla_{\eta'}\trans(\nabla_{\eta'}) 
\end{pmatrix} \\
&+\begin{pmatrix} 
0 & \frac{\gamma^2}{p'(\rho_p)}(\nu \Delta v_p+\tilde{\nu} \nabla  \mathrm{div} v_p)\\
0 & -\mathrm{div}v_p -\frac{1}{\rho_p}v_p\cdot \nabla_{\eta'}(\rho_p \cdot)+ \nabla v_p
\end{pmatrix}. 
\end{split}
\end{gather*}

Let $u^{(0)}$ and $u^{(0)*}$ denote the eigenfunctions for the eigenvalue $0$ of $-B_0$ and $-B_0^*$ satisfying 
\begin{equation*}
\langle\!\langle u^{(0)}, u^{(0)*}\rangle\!\rangle=1.
\end{equation*}
It then follows that 
\begin{equation*}
u_{\eta'}=\frac{1}{2 \pi i}\int_{|\lambda|=\frac{\pi}{4}} (\lambda + B_{\eta'})^{-1}u^{(0)} \,d\lambda
\end{equation*}
and 
\begin{equation*}
u_{\eta'}^*=\frac{1}{2 \pi i}\int_{|\lambda|=\frac{\pi}{4}} (\lambda + B_{\eta'}^*)^{-1}u^{(0)*} \,d\lambda 
\end{equation*}
are eigenfunctions for $B_{\eta'} \text{ and } B_{\eta'}^*$ associated with eigenvalue $\lambda_{\eta',0}$ and $\overline{\lambda_{\eta',0}}$, respectively.
Note that eigenfunctions for eigenvalues $\lambda_{\eta',k}$ are given by $e^{2\pi ikt}u_{\eta'}$ and the same holds for the adjoint eigenfunctions. 

We have the following estimates for the eigenfunctions for $u_{\eta'}$ and $u_{\eta'}^*$.
 \begin{thm} \label{u_eta}
Under the same assumptions of Theorem $\ref{critical eigen}$ the following estimates hold uniformly for $|\eta'| \leq r_0$ and $t\in\mathbb{T}_1:$
\begin{align*}
&\|u_{\eta'}(t)\|_{H^2(\Omega_{per})} \leq C,\\
&\no{u_{\eta'}(t)-u^{(0)}(t)}{}{H^2(\Omega_{per})}\leq C|\eta'|, \\
&\|u_{\eta'}^*(t)\|_{H^2(\Omega_{per})} \leq C.
\end{align*}
\end{thm}

Theorems \ref{critical eigen} and \ref{u_eta} will be proved in Section \ref{sec4}.

\section{Proof of Theorems \ref{critical eigen} and \ref{u_eta}}\label{sec4}

In this section, we prove Theorems \ref{critical eigen} and \ref{u_eta}. 
To do so, we consider the resolvent problem 
\begin{equation} \label{B eta}
(\lambda + B_{\eta'})u=F
\end{equation}
for $u \in D(B_{\eta'})$ with $|\eta'|\ll1$, where $F=\trans(f,g)$ is a given function. 

We expand $B_{\eta'}$ as 
\begin{equation*}
B_{\eta'}:= B_0+ \sum_{j=1}^{n-1}\eta_j B_j^{(1)} + \sum_{j,k=1}^{n-1}\eta_j \eta_k B_{j,k}^{(2)}.
\end{equation*}
Here
\begin{gather*}
\begin{split} 
B_0= \partial_t & +
\begin{pmatrix}
\mathrm{div}( v_p \cdot) & \gamma^2 \mathrm{div}( \rho_p \cdot)\\
\nabla \big(\frac{p'(\rho_p)}{\gamma^2 \rho_p}\cdot \big) & -\frac{\nu}{\rho_p}\Delta-\frac{\tilde{\nu}}{\rho_p}\nabla \mathrm{div}
\end{pmatrix} \\
&+\begin{pmatrix} 
0 & 0\\
\frac{1}{\gamma^2 \rho_p^2}(\nu\Delta v_p+\tilde{\nu}\nabla \mathrm{div} \  v_p) & v_p\cdot \nabla+\trans(\nabla v_p)
\end{pmatrix},
\end{split} \\
\begin{split}
&B_j^{(1)}=i
\begin{pmatrix}
v_p^j & \gamma^2 \rho_p \trans\mathbf{e}_j\\
 \big(\frac{p'(\rho_p)}{\gamma^2 \rho_p}\big)\mathbf{e}_j & -\frac{1}{\rho_p}(2 \nu \mathbf{e}_j \otimes \mathbf{e}_j \partial_{x_j} + \tilde{\nu}\mathbf{e}_j \mathrm{div} - \tilde{\nu} \nabla(\trans\mathbf{e}_j))+v_p^j
\end{pmatrix},
\end{split} \\
\begin{split} 
B_{j,k}^{(2)}=
\begin{pmatrix} 
0 & 0\\
0 & \frac{1}{\rho_p}(\nu \delta_{jk}I_n + \tilde{\nu}\mathbf{e}_j \trans\mathbf{e}_k)
\end{pmatrix}.
\end{split}
\end{gather*}
We set 
\begin{equation*}
M_{\eta'}=\sum_{j=1}^{n-1}\eta_j B_j^{(1)} + \sum_{j,k=1}^{n-1}\eta_j \eta_k B_{j,k}^{(2)}.
\end{equation*} 

We begin with investigating the spectral properties of $B_0$. 
For this purpose, we first consider the unique solvability for the time periodic problem 
\begin{gather}\label{prob L2}
\begin{cases}
\partial_tu+L(t)u=F, \\
\langle\phi(t)\rangle=0, 
\end{cases}
\end{gather}
when $F=\trans(f, g)\in L^2(\mathbb{T}_1; L^2_*(\Omega_{per})\times H^{-1}(\Omega_{per}))$. 

\begin{prop}\label{L2 time}
There exists positive constants $\nu_0$, $\gamma_0$, $\varepsilon_0$ and $a$ such that if $\frac{\nu^2}{\nu+\tilde{\nu}} \geq \nu_0$, $\frac{\gamma^2}{\nu+ \tilde{\nu}} \geq \gamma_0$ and $S\leq \varepsilon_{0}\frac{\nu^2}{\gamma^2\sqrt{\nu+\tilde{\nu}}}\sqrt{1-e^{-a\frac{\nu+\tilde{\nu}}{\gamma^2}}}$, then the following assertions hold true. 
For any $F\in L^2(\mathbb{T}_1; L^2_*(\Omega_{per})\times H^{-1}(\Omega_{per}))$, there exists a unique time periodic solution $u=\trans(\phi, w)\in C(\mathbb{T}_1; L_*^2(\Omega_{per})\times L^2(\Omega_{per}))\cap L^2(\mathbb{T}_1; L^2(\Omega_{per})\times H^1_0(\Omega_{per}))$ to $(\ref{prob L2})$. 
Furthermore, the solution $u$ satisfies 
\begin{equation}\label{est1 L2}
\begin{split}
&\|u(t)\|^2_{L^2(\Omega_{per}),\gamma} +\int^t_0 e^{-\beta_0(t-s)}\delta\|\phi(s)\|^2_{L^2(\Omega_{per})}\,ds \\
&+\int^t_0e^{-\beta_0(t-s)}\left(\nu\|\nabla w(s)\|^2_{L^2(\Omega_{per})}+\tilde{\nu}\|\mathrm{div}w(s)\|^2_{L^2(\Omega_{per})}\right)\,ds \\
&\leq \frac{C}{1-e^{-\beta_0}}\int^1_0\left(\frac{1}{\delta\gamma^4}+\frac{\delta^2}{\nu}\right)\|f(s)\|^2_{L^2(\Omega_{per})}+\frac{1}{\nu}\|g(s)\|^2_{H^{-1}(\Omega_{per})}\,ds. 
\end{split}
\end{equation}
Here $\delta = \frac{1}{4C}\min\left\{\frac{1}{\nu+\tilde{\nu}}, \frac{\nu}{\gamma^2}, \frac{1}{\gamma}\right\}$ and $\beta_0=C\min\left\{\nu+\tilde{\nu}, \delta\gamma^2\right\}$, where  $C$ is a positive constant independent of $\nu$, $\tilde{\nu}$, $\gamma$ and $S$.
\end{prop}

To prove Proposition \ref{L2 time}, we prepare the following lemma about the estimate of solution of the initial value problem for (\ref{prob L2}) under the initial condition 
\begin{equation}\label{IV L2}
u|_{t=0}=u_0=\trans(\phi_0, w_0), 
\end{equation}
when $u_0\in L^2_*(\Omega_{per})\times H^1_0(\Omega_{per})$ and $F\in L^2(\mathbb{T}_1; L^2_*(\Omega_{per})\times H^{-1}(\Omega_{per}))$. 
\begin{lem} \label{L2 global}
There exists positive constants $\nu_0$, $\gamma_0$, $\varepsilon_0$ and $a$ such that if $\frac{\nu^2}{\nu+\tilde{\nu}} \geq \nu_0$, $\frac{\gamma^2}{\nu+ \tilde{\nu}} \geq \gamma_0$ and $S\leq \varepsilon_{0}\frac{\nu^2}{\gamma^2\sqrt{\nu+\tilde{\nu}}}\sqrt{1-e^{-a\frac{\nu+\tilde{\nu}}{\gamma^2}}}$, then there exists a unique solution $u=\trans(\phi, w)\in C([0, \infty); L_*^2(\Omega_{per})\times L^2(\Omega_{per}))\cap L^2([0, \infty); L_*^2(\Omega_{per})\times H^1_0(\Omega_{per}))$ to $(\ref{prob L2})$ and $(\ref{IV L2})$. 
Furthermore, $u$ satisfies 
\begin{equation}\label{est2 L2}
\begin{split}
&\|u(t)\|^2_{L^2(\Omega_{per}),\gamma} +\int^t_0 e^{-\beta_0(t-s)}\delta\|\phi(s)\|^2_{L^2(\Omega_{per})}\,ds \\
&+\int^t_0e^{-\beta_0(t-s)}\left(\nu\|\nabla w(s)\|^2_{L^2(\Omega_{per})}+\tilde{\nu}\|\mathrm{div}w(s)\|^2_{L^2(\Omega_{per})}\right)\,ds \\
&\leq e^{-\beta_0t}\|u_0\|^2_{L^2(\Omega_{per}),\gamma} \\
&\hspace{4mm}+\frac{C}{1-e^{-\beta_0}}\int^1_0\left(\frac{1}{\delta\gamma^4}+\frac{\delta^2}{\nu}\right)\|f(s)\|^2_{L^2(\Omega_{per})}+\frac{1}{\nu}\|g(s)\|^2_{H^{-1}(\Omega_{per})}\,ds 
\end{split}
\end{equation}
for $t\geq0$. 
\end{lem}

\hspace{-6mm}{\bf Proof.} 
Since $v_p\in \cap^2_{j=0}C^j(\mathbb{T}_1; H^{4-2j}(\Omega_{per}))$, one can prove the existence of a solution to (\ref{prob L2}) and (\ref{IV L2}) with $u_0\in L^2_*(\Omega_{per})\times H^1_0(\Omega_{per})$ in a standard way by combining the method of characteristics and the parabolic theory. 

We prove the estimate (\ref{est2 L2}). 
We employ the energy method by Iooss-Padula \cite{Iooss-Padula}. 
We compute $\mathrm{Re}(\!(\partial_t u+L(t)u, u )\!)_t=\mathrm{Re}(\!(F, u)\!)_t$. 
In a similar way to the proof of \cite[Lemma 4.3]{E-K}, 
we see from Lemma \ref{lem Bogovskii} and Proposition \ref{timeper} that there exist positive constants $\nu_0$, $\gamma_0$, $\varepsilon_0$ and $a$ such that if $\delta=\frac{1}{4C}\min\{\frac{1}{\nu+\tilde{\nu}},\frac{\nu}{\gamma^2}, \frac{1}{\gamma}\}$, $\frac{\nu^2}{\nu+\tilde{\nu}}\geq \nu_0$, $\frac{\gamma^2}{\nu+\tilde{\nu}}\geq \gamma_0$ and $S\leq \varepsilon_{0}\frac{\nu^2}{\gamma^2\sqrt{\nu+\tilde{\nu}}}\sqrt{1-e^{-a\frac{\nu+\tilde{\nu}}{\gamma^2}}}$, then the following estimate holds: 
\begin{equation*}
\mathrm{Re}(\!(L(t)u,u)\!)_t\geq \nu\no{\nabla w}{2}{L^2(\Omega_{per})}+\tilde{\nu}\no{\mathrm{div}w}{2}{L^2(\Omega_{per})} +\delta\no{\phi}{2}{L^2(\Omega_{per})}.
\end{equation*}
On the other hand, we have 
\begin{align*}
\mathrm{Re}(\!(F,u)\!)_t\leq &\frac{\delta}{2}\no{\phi}{2}{L^2(\Omega_{per})}+\frac{\nu}{2}\no{\nabla w}{2}{L^2(\Omega_{per})} \\
&+C\left(\frac{1}{\delta\gamma^4}+\frac{\delta^2}{\nu}\right)\no{f}{2}{L^2(\Omega_{per})}+\frac{C}{\nu}\no{g}{2}{H^{-1}(\Omega_{per})}.
\end{align*}
It then follows 
\begin{align*}
&\frac{1}{2}\frac{d}{dt}(\!(u(t),u(t))\!)_t
+\frac{\nu}{2}\no{\nabla w}{2}{L^2(\Omega_{per})}+\tilde{\nu}\no{\mathrm{div}w}{2}{L^2(\Omega_{per})}+\frac{\delta}{2}\no{\phi}{2}{L^2(\Omega_{per})} \\
&\leq C\left(\frac{1}{\delta\gamma^4}+\frac{\delta^2}{\nu}\right)\no{f}{2}{L^2(\Omega_{per})}+\frac{C}{\nu}\no{g}{2}{H^{-1}(\Omega_{per})}.
\end{align*}
Multiplying this by $e^{\beta_0t}$ and integrating the resulting inequality over $[0,t]$, we have 
\begin{equation}\label{L2 decay}
\begin{split}
&\|u(t)\|^2_{L^2(\Omega_{per}),\gamma} +\int^t_0 e^{-\beta_0(t-s)}\delta\|\phi(s)\|^2_{L^2(\Omega_{per})}\,ds \\
&+\int^t_0e^{-\beta_0(t-s)}\left(\nu\|\nabla w(s)\|^2_{L^2(\Omega_{per})}+\tilde{\nu}\|\mathrm{div}w(s)\|^2_{L^2(\Omega_{per})}\right)\,ds \\
&\leq e^{-\beta_0t}\|u_0\|^2_{L^2(\Omega_{per}),\gamma}+C\int^t_0e^{-\beta_0(t-s)}\left(\frac{1}{\delta\gamma^4}+\frac{\delta^2}{\nu}\right)\|f(s)\|^2_{L^2(\Omega_{per})}\,ds \\
&\hspace{4mm}+C\int^t_0e^{-\beta_0(t-s)}\frac{1}{\nu}\|g(s)\|^2_{H^{-1}(\Omega_{per})}\,ds. 
\end{split}
\end{equation}
We apply Lemma \ref{intlemma} below to the right hand side of (\ref{L2 decay}) and obtain 
\begin{equation*}
\begin{split}
&\|u(t)\|^2_{L^2(\Omega_{per}),\gamma} +\int^t_0 e^{-\beta_0(t-s)}\delta\|\phi(s)\|^2_{L^2(\Omega_{per})}\,ds \\
&+\int^t_0e^{-\beta_0(t-s)}\left(\nu\|\nabla w(s)\|^2_{L^2(\Omega_{per})}+\tilde{\nu}\|\mathrm{div}w(s)\|^2_{L^2(\Omega_{per})}\right)\,ds \\
&\leq e^{-\beta_0t}\|u_0\|^2_{L^2(\Omega_{per}),\gamma}+\frac{C}{1-e^{-\beta_0}}\int^1_0\left(\frac{1}{\delta\gamma^4}+\frac{\delta^2}{\nu}\right)\|f(s)\|^2_{L^2(\Omega_{per})}\,ds \\
&\hspace{4mm}+\frac{C}{1-e^{-\beta_0}}\int^1_0\frac{1}{\nu}\|g(s)\|^2_{H^{-1}(\Omega_{per})}\,ds. 
\end{split}
\end{equation*}
This completes the proof. \qed 

\vspace{2ex}

\begin{lem} \label{intlemma}
If $f\in L^2(\mathbb{T}_1;L^2(\Omega_{per}))$, then 
\begin{equation*}
\int^t_0e^{-\beta_0(t-s)}\|f(s)\|^2_{L^2(\Omega_{per})}\,ds
\leq \frac{2}{1-e^{-\beta_0}}\int_0^1\|f(s)\|^2_{L^2(\Omega_{per})}\,ds. 
\end{equation*}
\end{lem}

\noindent 
{\bf Proof.} 
We set $N=[t]$, then 
\begin{align*}
&\int^t_0e^{-\beta_0(t-s)}\|f(s)\|^2_{L^2(\Omega_{per})}\,ds \\
&= \sum^{N-1}_{k=0}\int_k^{k+1}e^{-\beta_0(t-s)}\|f(s)\|^2_{L^2(\Omega_{per})}\,ds
+\int^t_N e^{-\beta_0(t-s)}\|f(s)\|^2_{L^2(\Omega_{per})}\,ds \\
&=:I_1 +I_2
\end{align*} 
Since $f\in L^2(\mathbb{T}_1;L^2(\Omega_{per}))$, we have
\begin{align*}
I_1
&= \sum^{N-1}_{k=0}\int_0^1e^{-\beta_0(t-\tilde{s}-k)}\|f(k+\tilde{s})\|^2_{L^2(\Omega_{per})}\,d\tilde{s} \\ 
&\leq \sum^{N-1}_{k=0}e^{-\beta_0(t-k-1)}\int_0^1\|f(\tilde{s})\|^2_{L^2(\Omega_{per})}\,d\tilde{s} \\
&\leq \frac{e^{\beta_0}}{e^{\beta_0}-1}\int_0^1\|f(\tilde{s})\|^2_{L^2(\Omega_{per})}\,d\tilde{s}
\end{align*}
and 
\begin{align*}
I_2
&=\int^{t-N}_0 e^{-\beta_0(t-\tilde{s}-N)}\|f(\tilde{s}+N)\|^2_{L^2(\Omega_{per})}\,d\tilde{s} \\ 
&\leq \int^{t-N}_0 e^{-\beta_0(t-\tilde{s}-N)}\|f(\tilde{s})\|^2_{L^2(\Omega_{per})}\,d\tilde{s} \\ 
&\leq \int^1_0 \|f(\tilde{s})\|^2_{L^2(\Omega_{per})}\,d\tilde{s}. 
\end{align*}
Therefore, we obtain 
\begin{align*}
\int^t_0e^{-\beta_0(t-s)}\|f(s)\|^2_{L^2(\Omega_{per})}\,ds 
&\leq \frac{2e^{\beta_0}}{e^{\beta_0}-1}\int_0^1\|f(\tilde{s})\|^2_{L^2(\Omega_{per})}\,d\tilde{s}. 
\end{align*}
This completes the proof. \qed

\vspace{2ex}

We are in a position to prove Proposition \ref{L2 time}. 

\vspace{2ex}

\hspace{-6mm}{\bf Proof of Proposition \ref{L2 time}.} 
We denote by $u^{\sharp}$ the solution of (\ref{prob L2}) and (\ref{IV L2}) with $u_0=0$. 
We then see from Lemma \ref{L2 global} that $u^{\sharp}$ satisfies 
\begin{equation} \label{est3 L2}
\begin{split}
&\|u^\sharp(t)\|^2_{L^2(\Omega_{per}),\gamma} \\
&\leq \frac{C}{1-e^{-\beta_0}}\int^1_0\left(\frac{1}{\delta\gamma^4}+\frac{\delta^2}{\nu}\right)\|f(s)\|^2_{L^2(\Omega_{per})}+\frac{1}{\nu}\|g(s)\|^2_{H^{-1}(\Omega_{per})}\,ds 
\end{split}
\end{equation}
 for $t\geq0$. 
Let $m,n \in\mathbb{N}$ with $m>n$. 
Since $F$ is periodic in $t$ of period $1$, the function $u^\sharp(t+(m-n))-u^\sharp(t)$ is the solution of (\ref{prob L2}) and (\ref{IV L2}) with $F=0$ and $u_0=u^{\sharp}(m-n)$. 
Hence, it follows from (\ref{est2 L2}) with $F=0$ that $u^\sharp(t+(m-n))-u^\sharp(t)$ satisfies 
\begin{equation*}
\|u^\sharp(t+(m-n))-u^\sharp(t)\|^2_{L^2(\Omega_{per}),\gamma}\leq e^{-\beta_0t}\|u^\sharp(m-n)\|^2_{L^2(\Omega_{per}),\gamma}. 
\end{equation*}
We set $t=n$ in this inequality. 
It then follows from (\ref{est3 L2}) that 
\begin{align*}
&\|u^\sharp(m)-u^\sharp(n)\|^2_{L^2(\Omega_{per}),\gamma} \\
&\leq e^{-\beta_0n}\|u^\sharp(m-n)\|^2_{L^2(\Omega_{per}),\gamma} \\
&\leq e^{-\beta_0n}\frac{C}{1-e^{-\beta_0}}\int^1_0\left(\frac{1}{\delta\gamma^4}+\frac{\delta^2}{\nu}\right)\|f(s)\|^2_{L^2(\Omega_{per})}+\frac{1}{\nu}\|g(s)\|^2_{H^{-1}(\Omega_{per})}\,ds. 
\end{align*}
We thus obtain 
\begin{equation*}
\|u^\sharp(m)-u^\sharp(n)\|^2_{L^2(\Omega_{per}),\gamma} \rightarrow 0 \quad (n\rightarrow \infty).
\end{equation*}
Therefore, $\{u^\sharp(m)\}$ is a Cauchy sequence in $L_*^2(\Omega_{per})\times L^2(\Omega_{per})$. 
It then follows that there exists $\tilde{u}^\sharp \in L_*^2(\Omega_{per})\times L^2(\Omega_{per})$ such that 
$u^\sharp(m)$ converges to $\tilde{u}^\sharp$ strongly in $L_*^2(\Omega_{per})\times L^2(\Omega_{per})$, and $\tilde{u}^\sharp$ satisfies   
\begin{equation} \label{est4 L2}
\|\tilde{u}^\sharp\|^2_{L^2(\Omega_{per}),\gamma}\leq \frac{C}{1-e^{-\beta_0}}\int^1_0\left(\frac{1}{\delta\gamma^4}+\frac{\delta^2}{\nu}\right)\|f(s)\|^2_{L^2(\Omega_{per})}+\frac{1}{\nu}\|g(s)\|^2_{H^{-1}(\Omega_{per})}\,ds. 
\end{equation}
We then see from the argument by Valli that the solution $u$ of (\ref{prob L2}) and (\ref{IV L2}) with $u_0=\tilde{u}^{\sharp}$ is a time periodic solution of (\ref{prob L2}).
Furthermore, applying (\ref{est2 L2}) and (\ref{est4 L2}), we obtain  
\begin{equation*} 
\begin{split}
&\|u(t)\|^2_{L^2(\Omega_{per}),\gamma} +\int^t_0 e^{-\beta_0(t-s)}\delta\|\phi(s)\|^2_{L^2(\Omega_{per})}\,ds \\
&+\int^t_0e^{-\beta_0(t-s)}\left(\nu\|\nabla w(s)\|^2_{L^2(\Omega_{per})}+\tilde{\nu}\|\mathrm{div}w(s)\|^2_{L^2(\Omega_{per})}\right)\,ds \\
&\leq \frac{C}{1-e^{-\beta_0}}\int^1_0\left(\frac{1}{\delta\gamma^4}+\frac{\delta^2}{\nu}\right)\|f(s)\|^2_{L^2(\Omega_{per})}+\frac{1}{\nu}\|g(s)\|^2_{H^{-1}(\Omega_{per})}\,ds. 
\end{split}
\end{equation*}
This completes the proof. \qed 

\vspace{2ex}

The following proposition shows that $0$ is an eigenvalue of $-B_0$. We also give the estimates of an eigenfunction for the eigenvalue 0.

 \begin{prop} \label{exist}
There exist positive constants $\nu_0$, $\gamma_0$, $\varepsilon_0$ and $a$ such that the following assertions hold. 
If $\frac{\nu^2}{\nu+\tilde{\nu}} \geq \nu_0$, $\frac{\gamma^2}{\nu+ \tilde{\nu}} \geq \gamma_0$ and $S\leq \varepsilon_{0}\frac{\nu^2}{\gamma^2\sqrt{\nu+\tilde{\nu}}}\sqrt{1-e^{-a\frac{\nu+\tilde{\nu}}{\gamma^2}}}$, then there exists a solution $u^{(0)}=\trans(\phi^{(0)},w^{(0)})\in D(B_0)$ of 
\begin{equation}
\begin{cases}
B_0u^{(0)}=0,\\
\mean{\phi^{(0)}}=1.
\end{cases}\label{4.17}
\end{equation} 
Furthermore, $\tilde{u}^{(0)}=\trans(\tilde{\phi}^{(0)}, w^{(0)})=(\phi^{(0)}-1,w^{(0)})$ satisfies
 \begin{gather}
\frac{1}{\gamma^2}\nom{\tilde{\phi}^{(0)}(t)}{2}{2} +\nom{w^{(0)}(t)}{2}{2} \leq \frac{C}{\nu\gamma^2}, \label{4.18a}\\
 \int_{0}^{1} \frac{1}{\nu+\tilde{\nu}}\|\tilde{\phi}^{(0)}\|_{H^2(\Omega_{per})}^2+\frac{\nu+\tilde{\nu}}{\gamma^4}\|\partial_t\tilde{\phi}^{(0)}\|_{L^2(\Omega_{per})}^2 + \frac{\nu^2}{\nu+\tilde{\nu}} \nom{w^{(0)}}{2}{3} \,ds 
 \leq \frac{C}{\nu\gamma^2}.\label{4.18}
 \end{gather} 
 \end{prop}

Before proving Proposition \ref{exist}, we state one proposition on the spectrum of $-B_0$ which immediately follows from Proposition \ref{exist}.
\begin{prop}
Under the assumption of Proposition $\ref{exist}$ for each $k\in\mathbb{Z}$, $2\pi ik$ is an eigenvalue of $-B_0$ with eigenfunction $e^{2\pi ikt}u^{(0)}$. 
\end{prop}

To prove Proposition \ref{exist}, we decompose $\phi$ into $\phi=1+\tilde{\phi}$ and rewrite (\ref{4.17}) for $u=\trans(\phi, w)$ as (\ref{prob L2}) for $\tilde{u}=\trans(\tilde{\phi}, w)$. 
We thus consider the time periodic problem for (\ref{prob L2}) with 
\begin{equation}\label{eigen F}
F=
\begin{pmatrix}
-\mathrm{div}v_p \\
\nabla\left(\frac{P'(\rho_p)}{\gamma^2\rho_p} \right)+\frac{1}{\gamma^2\rho_p^2}(\nu\Delta v_p+\tilde{\nu}\nabla\mathrm{div}v_p)
\end{pmatrix}.
\end{equation}
As in the proof of Proposition \ref{L2 time}, we first consider the initial value problem for (\ref{prob L2}) with $F$ given in (\ref{eigen F}) under the initial condition (\ref{IV L2}) with $u_0\in H^2_*(\Omega_{per})\times (H^2(\Omega_{per})\cap H^1_0(\Omega_{per}))$.

\begin{lem}\label{H2 global}
There exist positive constants $\nu_0$, $\gamma_0$, $\varepsilon_0$, $a$ and $C_E$ such that if $\frac{\nu^2}{\nu+\tilde{\nu}} \geq \nu_0$, $\frac{\gamma^2}{\nu+ \tilde{\nu}} \geq \gamma_0$ and $S\leq \varepsilon_{0}\frac{\nu^2}{\gamma^2\sqrt{\nu+\tilde{\nu}}}\sqrt{1-e^{-a\frac{\nu+\tilde{\nu}}{\gamma^2}}}$, then there exists a unique solution $u=\trans(\phi, w)\in \cap_{j=0}^1C^j([0, \infty); H_*^{2-2j}(\Omega_{per})\times H^{2-2j}(\Omega_{per}))\cap H^j([0, \infty); H_*^{2-2j}(\Omega_{per})\times (H^{3-2j}(\Omega_{per})\cap H^1_0(\Omega_{per})))$ to $(\ref{prob L2})$ and $(\ref{IV L2})$ with $u_0\in H^2_*(\Omega_{per})\times (H^2(\Omega_{per})\cap H^1_0(\Omega_{per}))$ and $F$ given by $(\ref{eigen F})$. 
Furthermore, $u$ satisfies 
\begin{align*}
&\frac{1}{\gamma^2}\nom{\phi(t)}{2}{2}+\nom{w(t)}{2}{2}+\int_0^te^{-a\frac{\nu+\tilde{\nu}}{\gamma^2}(t-s)}\frac{\nu^2}{\nu+\tilde{\nu}}\nom{w(s)}{2}{3}\,ds \\
&+\int_0^te^{-a\frac{\nu+\tilde{\nu}}{\gamma^2}(t-s)}\left(\frac{1}{\nu+\tilde{\nu}}\|\phi(s)\|^2_{H^2(\Omega_{per})}+\frac{\nu+\tilde{\nu}}{\gamma^4}\|\partial_t\phi(s)\|^2_{L^2(\Omega_{per})}\right)\,ds \\
&\leq (1+C_E) \left\{ e^{-a\frac{\nu+\tilde{\nu}}{\gamma^2}t}\left(\frac{1}{\gamma^2}\nom{\phi_0}{2}{2}+\nom{w_0}{2}{2}\right) 
+\frac{C}{1-e^{-a\frac{\nu+\tilde{\nu}}{\gamma^2}}}\frac{1}{\gamma^4} \right\}. 
\end{align*}
\end{lem}

\hspace{-6mm}{\bf Proof.} 
Since $v_p\in \cap_{j=0}^2C^j(\mathbb{T}_1; H^{4-2j}(\Omega_{per}))$, one can prove the existence of solution $u$ to (\ref{prob L2}) and (\ref{IV L2}) with $u_0\in H^2_*(\Omega_{per})\times (H^2(\Omega_{per})\cap H^1_0(\Omega_{per}))$ 
in a standard way by combining the method of characteristics and the parabolic theory. The estimate of the solution $u$ is obtained in a similar manner to the proof of Proposition \ref{apriori} below. We here give an outline of the proof of the estimate. 
 
We rewrite $L(t)$ as  
\begin{equation*}
L(t)=A+M(t),
\end{equation*}
where 
\begin{gather*}
A=
\begin{pmatrix}
0 & \gamma^2\mathrm{div} \\
\nabla & -\nu\Delta-\tilde{\nu}\nabla\mathrm{div}
\end{pmatrix}, \\
\begin{split}
M(t)=&
\begin{pmatrix}
v_p(t)\cdot\nabla+\mathrm{div}v_p(t) & \gamma^2\mathrm{div}(\phi_p(t)\cdot) \\
\nabla(p^{(1)}(\phi_p(t))\phi_p(t)\cdot) & \frac{\phi_p(t)}{1+\phi_p(t)}(\nu\Delta+\tilde{\nu}\nabla\mathrm{div})
\end{pmatrix} \\
&+\begin{pmatrix}
0 & 0 \\
\frac{1}{\gamma^2\rho_p^2(t)}(\nu\Delta v_p(t)+\tilde{\nu}\nabla\mathrm{div}v_p(t)) & v_p(t)\cdot\nabla+\trans(\nabla v_p(t))
\end{pmatrix}.
\end{split}
\end{gather*}
Here $p^{(1)}(\phi)=\int^1_0\frac{p''(1+\theta\phi)}{\gamma^2(1+\theta\phi)}d\theta$. 
Furthermore, we rewrite (\ref{prob L2}) as 
\begin{equation*}
\partial_tu+Au=F-Mu.
\end{equation*}
As in the proof of Proposition \ref{energyest}, applying Proposition \ref{basicest} with $m=2$ for 
\begin{equation}\label{eigen F1}
\begin{pmatrix} 
f \\
g
\end{pmatrix} 
=F-Mu, \qquad 
\tilde{f}=\phi\mathrm{div} v_p +\gamma^2\mathrm{div}(\phi_p w). 
\end{equation}
one can show that there exists a positive constant $\nu_0$ such that if $\frac{\nu^2}{\nu+\tilde{\nu}}\geq\nu_0$, then 
\begin{gather*}
\frac{d}{dt}\tilde{E}_2(u)+D_2(u)\leq \tilde{N}_2(u). 
\end{gather*}
Here $\tilde{E}_{2}(u)$ and $D_{2}(u)$ are the same functionals as those given in Propositions  \ref{energyest} and \ref{est}, respectively; and $\tilde{N}_{2}$ is a functional satisfying 
\begin{align*}
\tilde{N}_2(u)\leq C\frac{\gamma^2}{\nu+\tilde{\nu}}[\![v_p]\!]_4D_2(u)+\frac{1}{4}D_2(u) 
+C\frac{1}{\nu+\tilde{\nu}}\nom{v_p}{2}{5}+C\frac{1}{\nu}\nom{\nabla\phi_p}{2}{3}. 
\end{align*}
It then follows that there exist positive constants $\nu_0$, $\gamma_0$, $\varepsilon_0$ and $a$ such that if $\frac{\nu^2}{\nu+\tilde{\nu}}\geq \nu_0$, $\frac{\gamma^2}{\nu+\tilde{\nu}}\geq \gamma_0$ and $S\leq \varepsilon_{0}\frac{\nu^2}{\gamma^2\sqrt{\nu+\tilde{\nu}}}\sqrt{1-e^{-a\frac{\nu+\tilde{\nu}}{\gamma^2}}}$, we have  
\begin{equation}\label{4.20}
\begin{split}
\frac{d}{dt}\tilde{E}_2(u)+D_2(u)
\leq &C\frac{1}{\nu+\tilde{\nu}}\nom{v_p}{2}{5}+C\frac{1}{\nu}\nom{\nabla\phi_p}{2}{3}.
\end{split} 
\end{equation}
Applying now the argument of the proof of Proposition \ref{apriori} below, we have 
\begin{align*}
&E_2(u(t))+ \int_0^te^{-a\frac{\nu+\tilde{\nu}}{\gamma^2}(t-s)}D_2(u(s))\,ds  \\
&\leq (1+C_E) \left\{ e^{-a\frac{\nu+\tilde{\nu}}{\gamma^2}t}E_2(u_0) 
+C\int^1_0e^{-a\frac{\nu+\tilde{\nu}}{\gamma^2}(t-s)}\left(\frac{1}{\nu+\tilde{\nu}}\nom{v_p}{2}{5}+\frac{1}{\nu}\nom{\nabla\phi_p}{2}{3}\right)\,ds \right\}. 
\end{align*}
Here $E_{2}(u)$ is the same functional given in Proposition \ref{est} below; and $C_E$ is the positive constant given in Lemma \ref{extension} below. 
This, together with Lemma \ref{intlemma}, implies that 
\begin{equation*}
\begin{split}
&E_2(u(t))+ \int_0^te^{-a\frac{\nu+\tilde{\nu}}{\gamma^2}(t-s)}D_2(u(s))\,ds  \\
&\leq (1+C_E) \left\{ e^{-a\frac{\nu+\tilde{\nu}}{\gamma^2}t}E_2(u_0) 
+\frac{C}{1-e^{-a\frac{\nu+\tilde{\nu}}{\gamma^2}}}\int^1_0\left(\frac{1}{\nu+\tilde{\nu}}\nom{v_p}{2}{5}+\frac{1}{\nu}\nom{\nabla\phi_p}{2}{3}\right)\,ds \right\}. 
\end{split}
\end{equation*}
The desired estimate follows from this inequality by applying Proposition \ref{timeper} to the second term of the right-hand side. This completes the proof. \qed 

\vspace{2ex}

We are in a position to prove Proposition \ref{exist}. 

\vspace{2ex}

\hspace{-6mm}{\bf Proof of Proposition \ref{exist}.}
Decomposing $\phi$ into $\phi=1+\tilde{\phi}$, we rewrite the problem (\ref{4.17}) for $u=\trans(\phi, w)$ as the problem (\ref{prob L2}) with $F$ given in (\ref{eigen F}) for $\tilde{u}=\trans(\tilde{\phi}, w)$.  
Based on Lemma \ref{H2 global}, in a similar manner to the proof of Proposition \ref{L2 time}, by using the argument of Valli \cite{Valli}, we can obtain a time periodic solution $\tilde{u}=\trans(\tilde{\phi}, w)$ to (\ref{prob L2}) with $F$ given in (\ref{eigen F}) satisfying 
\begin{align*}
\frac{1}{\gamma^2}\nom{\tilde{\phi}(t)}{2}{2}+\nom{w(t)}{2}{2} 
\leq \frac{C}{1-e^{-a\frac{\nu+\tilde{\nu}}{\gamma^2}}}\frac{1}{\gamma^4}
\end{align*}
and 
\begin{align*}
\int_0^1\left(\frac{1}{\nu+\tilde{\nu}}\|\tilde{\phi}\|^2_{H^2(\Omega_{per})}+\frac{\nu+\tilde{\nu}}{\gamma^4}\|\partial_t\tilde{\phi}\|^2_{L^2(\Omega_{per})}+\frac{\nu^2}{\nu+\tilde{\nu}}\nom{w}{2}{3} \right)\,ds
\leq \frac{C}{1-e^{-a\frac{\nu+\tilde{\nu}}{\gamma^2}}}\frac{1}{\gamma^4}.
\end{align*}
Since $1-e^{-a\frac{\nu+\tilde{\nu}}{\gamma^2}}\geq \frac{a}{2}\frac{\nu+\tilde{\nu}}{\gamma^2}$ if $\frac{\gamma^2}{\nu+\tilde{\nu}}\geq\gamma_0$ for some positive constant $\gamma_0$, we have 
\begin{align*}
\frac{1}{\gamma^2}\nom{\tilde{\phi}(t)}{2}{2}+\nom{w(t)}{2}{2} 
\leq \frac{C}{\nu\gamma^2}
\end{align*}
and 
\begin{align*}
\int_0^1\left(\frac{1}{\nu+\tilde{\nu}}\|\tilde{\phi}\|^2_{H^2(\Omega_{per})}+\frac{\nu+\tilde{\nu}}{\gamma^4}\|\partial_t\tilde{\phi}\|^2_{L^2(\Omega_{per})}+\frac{\nu^2}{\nu+\tilde{\nu}}\nom{w}{2}{3} \right)\,ds
\leq \frac{C}{\nu\gamma^2}.
\end{align*}
This completes the proof. \qed 

\vspace{2ex}

To prove that $0$ is a simple eigenvalue of $-B_0$, we prepare the following lemma. 
  \begin{lem}\label{lem u*}
  Let $\Pi^{(0)}$ be defined by  
$$\Pi^{(0)}u=\inn{u}{u^{(0)*}}u^{(0)}=\mean{\phi}u^{(0)}$$ 
for $u=\trans(\phi, w)$, where 
  $$u^{(0)*}=\gamma^2 \begin{pmatrix}
  \frac{\gamma^2 \rho_p}{P'(\rho_p)}\\ 0
  \end{pmatrix}.$$
  Then the following assertions hold.
  \begin{enumerate} [label=(\roman*)]
  \item[{\rm (i)}] $u^{(0)*}$ satisfies $u^{(0)*} \in D(B_0^*)$, $B^*_0u^{(0)*}=0$ and $\langle u^{(0)}(t),u^{(0)*}(t)\rangle =1$.
  \item[{\rm (ii)}] $\Pi^{(0)}$ is a bounded projection on $X$ satisfying 
\begin{equation*}
\Pi^{(0)}B_0\subset B_0\Pi^{(0)}=0, \ 
\Pi^{(0)}u^{(0)}=u^{(0)} \text{ and }\ \Pi^{(0)}X=\mathrm{span}\{u^{(0)} \}.
\end{equation*}
\end{enumerate}
\end{lem}
One can prove Lemma \ref{lem u*} by straightforward computations. 

\vspace{2ex}

We now prove the simplicity of the eigenvalue of $-B_0$. 
Let $X_0$ and $X_1$ be defined by 
\begin{equation*}
X_0=\Pi^{(0)}X~ \text{and}\ X_1=(I-\Pi^{(0)})X.
\end{equation*}
Observe that $u=\trans(\phi,w)\in X_1$ if and only if $\mean{\phi}=0$.
 
As for $X_0$ and $X_1$, it holds the following assertions.
 
 \begin{prop} \label{eigenspace}
If $\frac{\nu^2}{\nu+\tilde{\nu}} \geq \nu_0$, $\frac{\gamma^2}{\nu+ \tilde{\nu}} \geq \gamma_0$ and $S\leq \varepsilon_{0}\frac{\nu^2}{\gamma^2\sqrt{\nu+\tilde{\nu}}}\sqrt{1-e^{-a\frac{\nu+\tilde{\nu}}{\gamma^2}}}$, then 
 \begin{enumerate}[label=(\roman*)]
 \item[{\rm (i)}] $X_0=\mathrm{Ker}(B_0) \text{ and }X_1= \mathrm{Ran}(B_0);$ $X_1$ is closed.
 \item[{\rm (ii)}] $X=X_0\oplus X_1$.
 \item[{\rm (iii)}] 0 is a simple eigenvalue of $-B_0.$
 \end{enumerate}
 \end{prop}

  \hspace{-6mm}{\bf Proof.}
Let us show $X_1=\mathrm{Ran}(B_0)$. 
We first assume that $F=\trans(f,g) \in \mathrm{Ran}(B_0)$. 
There exists a function $u=\trans(\phi,w) \in D(B_0)$ such that $B_0 u=F$. Applying $\Pi^{(0)}$ to $B_0 u=F$, we have $$ 0=\Pi^{(0)}B_0 u=\Pi^{(0)}F=\mean{f} u^{(0)}.$$
    This implies that $\mean{f}=0$ and hence $F \in X_1$. 
We thus obtain $\mathrm{Ran}(B_0) \subset X_1$. 

We next prove $X_1\subset\mathrm{Ran}(B_0)$. 
Let $F=\trans(f,g) \in X_1$. 
We will show that there exists a unique solution $u\in D(B_0)\cap X_1$ to 
\begin{equation}
B_0u=F. \label{aa}
\end{equation} 
We define $\tilde{\Pi}^{(0)}$ by $\tilde{\Pi}^{(0)}u=\trans(\langle\phi\rangle, 0)$ for $u=\trans(\phi, w)$. 
We decompose $u=\trans(\phi,w)$ as 
\begin{align*}
u=\tilde{\Pi}^{(0)}u + (I-\tilde{\Pi}^{(0)})u=\trans(\langle\phi\rangle, 0) + u_1, 
\end{align*}
where $u_1=\trans(\phi_1, w_1)\in L^2(\mathbb{T}_1; L^2_*(\Omega_{per})\times L^2(\Omega_{per}))$. 
Applying $\tilde{\Pi}^{(0)}$ and $(I-\tilde{\Pi}^{(0)})$ to $B_0u=F$, we have 
\begin{gather} \label{<phi0>}
\partial_t \langle \phi \rangle = \langle f \rangle, \\ \label{u_10}
\partial_t u_1 + L(t)u_1 = F_1  -L(t)\trans(\langle \phi \rangle , 0)
\end{gather}
since $\tilde{\Pi}^{(0)}L(t)u=0$,
    where $F_1=\trans(f_1,g_1):=(I-\tilde{\Pi}^{(0)})F$. 
Integrating (\ref{<phi0>}) in $[0,t]$, we have 
\begin{equation*}
\langle \phi(t)\rangle=\langle\phi(0)\rangle +\int^t_0\langle f(s)\rangle\,ds. 
\end{equation*}
Furthermore, we determine $\langle\phi(0)\rangle=\int^1_0s\langle f(s)\rangle\,ds$ so that $\int^1_0\langle\phi(t)\rangle\,dt=0$. 
Consequently, we obtain 
\begin{equation}
\langle\phi(t)\rangle = \int^1_0s\langle f(s)\rangle\,ds+\int^t_0\langle f(s)\rangle\,ds. \label{<phi1>}
\end{equation}
As for (\ref{u_10}), it follows from Proposition \ref{L2 time} that there exist positive constants $\nu_0$, $\gamma_0$, $\varepsilon_0$ and $a$ such that if $\delta=\frac{1}{4C}\min\{\frac{1}{\nu+\tilde{\nu}},\frac{\nu}{\gamma^2}, \frac{1}{\gamma}\}$, $\frac{\nu^2}{\nu+\tilde{\nu}}\geq \nu_0$, $\frac{\gamma^2}{\nu+\tilde{\nu}}\geq \gamma_0$ and $S\leq \varepsilon_{0}\frac{\nu^2}{\gamma^2\sqrt{\nu+\tilde{\nu}}}\sqrt{1-e^{-a\frac{\nu+\tilde{\nu}}{\gamma^2}}}$, there exists a time periodic solution $u_1=\trans(\phi_1, w_1)\in L^2(\mathbb{T}_1; L^2_*(\Omega_{per})\times H^1_0(\Omega_{per}))$ to (\ref{u_10}) satisfying 
\begin{equation}\label{est-u1}
\begin{split}
&\|u(t)\|^2_{L^2(\Omega_{per}),\gamma} +\int^t_0 e^{-\beta_0(t-s)}\delta\|\phi(s)\|^2_{L^2(\Omega_{per})}\,ds \\
&+\int^t_0e^{-\beta_0(t-s)}\left(\nu\|\nabla w(s)\|^2_{L^2(\Omega_{per})}+\tilde{\nu}\|\mathrm{div}w(s)\|^2_{L^2(\Omega_{per})}\right)\,ds \\
&\leq \frac{C}{1-e^{-\beta_0}}\int^1_0\left(\frac{1}{\delta\gamma^4}+\frac{\delta^2}{\nu}\right)\|f(s)\|^2_{L^2(\Omega_{per})}+\frac{1}{\nu}\|g(s)\|^2_{H^{-1}(\Omega_{per})}\,ds. 
\end{split}
\end{equation}
Hence, there exists a unique solution $u=\trans(\phi,w)\in D(B_0)\cap X_1$ to (\ref{aa}). 
This shows $X_1 \subset \mathrm{Ran}(B_0)$. Therefore $X_1=\mathrm{Ran}(B_0)$.

Let us show $X_0=\mathrm{Ker}(B_0)$. 
We assume that $u$ is the solution of $B_0u=0$. 
Decomposing $u$ into $u=u_0+u_1$ with $u_j\in X_j$ for $j=0,1$, we have $B_0u_j=0$. 
It follows from the previous argument that $u_1$ is a unique solution to $B_0u_1=0$, 
and  we see from (\ref{est-u1}) with $F=0$ that $u_1=0$.  
Consequently, it holds that $u=u_0\in X_0$ and $\mathrm{Ker}(B_0)\subset X_0$. Therefore, $\mathrm{Ker}(B_0)=X_0$.
This completes the proof.
 \qed

\begin{rem}\label{rem simple}
One can show that, for each $k\in\mathbb{Z}$, $2\pi ik$ is a simple eigenvalue of $-B_0$. 
\end{rem}

We next establish the resolvent estimate for $\eta' =0$. 
We consider 
\begin{equation} \label{range}
 \lambda u+ \partial_tu +L(t)u=F, 
\end{equation}
where $F=\trans(f,g)\in L^2(\mathbb{T}_1; L^2(\Omega_{per})\times H^{-1}(\Omega_{per}))$.

\begin{prop} \label{prop C}
There exist positive constants $\nu_0$, $\gamma_0$, $\varepsilon_0$ and $a$ such that if $\frac{\nu^2}{\nu+\tilde{\nu}} \geq \nu_0$, $\frac{\gamma^2}{\nu+ \tilde{\nu}} \geq \gamma_0$ and $S\leq \varepsilon_{0}\frac{\nu^2}{\gamma^2\sqrt{\nu+\tilde{\nu}}}\sqrt{1-e^{-a\frac{\nu+\tilde{\nu}}{\gamma^2}}}$, then 
there exists a unique time periodic solution $u=\trans(\phi, w)$ of $(\ref{range})$ for $\lambda\in\mathbb{C}$ with $\mathrm{Re}\lambda > -\frac{\beta_0}{2}$ and $\lambda\neq2\pi ik$ $(k\in\mathbb{Z})$. 
Furthermore, $u$ satisfies 
\begin{equation*}
\begin{split}
&\|u(t)\|^2_{L^2(\Omega_{per}),\gamma}+\int_0^te^{-(2\mathrm{Re}\lambda+\beta_0)(t-s)}\|(\phi-\langle\phi\rangle)(s)\|^2_{L^2(\Omega_{per})}\,ds \\
&+\int_0^te^{-(2\mathrm{Re}\lambda+\beta_0)(t-s)}\left(\|\nabla w(s)\|^2_{L^2(\Omega_{per})}+\|\mathrm{div}w(s)\|^2_{L^2(\Omega_{per})}\right)\,ds  \\
&\leq \frac{C}{(2\mathrm{Re}\lambda+\beta_0)|1-e^{-\lambda}|^2}\int^1_0\|f(s)\|^2_{L^2(\Omega_{per})}\,ds \\
&\hspace{4mm}+\frac{C}{1-e^{-(2\mathrm{Re}\lambda+\beta_0)}}\int^1_0\|F(s)\|^2_{L^2(\Omega_{per})\times H^{-1}(\Omega_{per})}\,ds 
\end{split}
\end{equation*}
for $t\in \mathbb{T}_1$.
\end{prop}
\hspace{-6mm}{\bf Proof.}
As in the proof of Proposition \ref{eigenspace}, we apply $\tilde{\Pi}^{(0)}$ and $(I-\tilde{\Pi}^{(0)})$ to \eqref{range}. Then we have,
\begin{equation} \label{<phi>}
\lambda \langle \phi \rangle+\partial_t \langle \phi \rangle  = \langle f \rangle 
\end{equation}
\begin{equation} \label{u_1}
\lambda u_1+\partial_t u_1 + L(t)u_1 = G_1,
\end{equation}
where $G_1=-L(t)\trans(\langle\phi\rangle,0)+F_1$. 
 Since we look for a time periodic solution, $\langle \phi \rangle$ must satisfy  
\begin{equation} \label{phi_1}
\begin{cases}
\langle \phi \rangle (t)= e^{-\lambda t}\langle \phi \rangle (0)+ \int_{0}^{t}e^{-\lambda (t-s)}\langle f(s) \rangle\,ds, \\
\langle \phi \rangle (1)= \langle \phi \rangle (0). 
\end{cases}
\end{equation}
In (\ref{phi_1}), set $t=1$. 
We then obtain 
\begin{equation*}
\langle \phi \rangle (0)= e^{-\lambda}\langle \phi \rangle (0)+ \int_{0}^{1}e^{-\lambda (1-s)}\langle f(s) \rangle\, ds.
\end{equation*}
Therefore, if $1-e^{-\lambda} \neq 0$, namely, if $\lambda \neq 2 \pi ik$ for $k \in \mathbb{Z}$, then
\begin{equation}\label{phi3}
\langle \phi \rangle (0)=\frac{1}{1-e^{-\lambda}}  \int_{0}^{1}e^{-\lambda (1-s)}\langle f(s) \rangle\, ds.
\end{equation}
Substituting (\ref{phi3}) to the first equation of (\ref{phi_1}), we have 
\begin{equation*}
\langle \phi \rangle (t)=\frac{e^{-\lambda t}}{1-e^{-\lambda}}  \int_{0}^{1}e^{-\lambda (1-s)}\langle f(s) \rangle\, ds + \int_{0}^{t}e^{-\lambda (t-s)}\langle f(s) \rangle\, ds.
\end{equation*} 
Therefore, if $\lambda\in\Sigma$, then we obtain 
\begin{equation}\label{phi2}
|\langle \phi \rangle (t)| \leq \frac{C}{|1-e^{-\lambda}|} \int_{0}^{1} \|f\|_{L^2(\Omega_{per})}\, dt.
\end{equation}

We next consider (\ref{u_1}). 
We set $v(t)=e^{\lambda t}u(t)$. 
Since $v$ satisfies (\ref{prob L2}) with $F=e^{\lambda t}G_1$ and the estimate (\ref{est2 L2}), there exist positive constants $\nu_0$, $\gamma_0$, $\varepsilon_0$ and $a$ such that if $\frac{\nu^2}{\nu+\tilde{\nu}}\geq \nu_0$, $\frac{\gamma^2}{\nu+\tilde{\nu}}\geq \gamma_0$ and $S\leq \varepsilon_{0}\frac{\nu^2}{\gamma^2\sqrt{\nu+\tilde{\nu}}}\sqrt{1-e^{-a\frac{\nu+\tilde{\nu}}{\gamma^2}}}$, we have  
\begin{equation*}
\begin{split}
&\|u_1(t)\|^2_{L^2(\Omega_{per}),\gamma}+\int_0^te^{-(2\mathrm{Re}\lambda+\beta_0)(t-s)}\|\phi_1(s)\|^2_{L^2(\Omega_{per})}\,ds \\
&+\int_0^te^{-(2\mathrm{Re}\lambda+\beta_0)(t-s)}\left(\|\nabla w_1(s)\|^2_{L^2(\Omega_{per})}+\|\mathrm{div}w_1(s)\|^2_{L^2(\Omega_{per})}\right)\,ds  \\
&\leq e^{-(2\mathrm{Re}\lambda+\beta_0)t}\|u_1(0)\|^2_{L^2(\Omega_{per})} \\
&\hspace{4mm}+C\int^t_0e^{-(2\mathrm{Re}\lambda+\beta_0)(t-s)}\|G_1(s)\|^2_{L^2(\Omega_{per})\times H^{-1}(\Omega_{per})}\,ds. 
\end{split}
\end{equation*}
As in the proofs of Proposition \ref{L2 time} and Lemma \ref{L2 global}, one can see that if $\mathrm{Re}\lambda + \frac{\beta_0}{2} > 0$, then there exists a time periodic solution $u_1=\trans(\phi_1, w_1)$ to (\ref{u_1}) satisfying 
\begin{equation*}
\begin{split}
&\|u_1(t)\|^2_{L^2(\Omega_{per}),\gamma}+\int_0^te^{-(2\mathrm{Re}\lambda+\beta_0)(t-s)}\|\phi_1(s)\|^2_{L^2(\Omega_{per})}\,ds \\
&+\int_0^te^{-(2\mathrm{Re}\lambda+\beta_0)(t-s)}\left(\|\nabla w_1(s)\|^2_{L^2(\Omega_{per})}+\|\mathrm{div}w_1(s)\|^2_{L^2(\Omega_{per})}\right)\,ds  \\
&\leq C\int^t_0e^{-(2\mathrm{Re}\lambda+\beta_0)(t-s)}\|G_1(s)\|^2_{L^2(\Omega_{per})\times H^{-1}(\Omega_{per})}\,ds \\
&\leq \frac{C}{(2\mathrm{Re}\lambda+\beta_0)|1-e^{-\lambda}|^2}\int^1_0\|f(s)\|^2_{L^2(\Omega_{per})}\,ds \\
&\hspace{4mm}+C\int^t_0e^{-(2\mathrm{Re}\lambda+\beta_0)(t-s)}\|F_1(s)\|^2_{L^2(\Omega_{per})\times H^{-1}(\Omega_{per})}\,ds. 
\end{split}
\end{equation*}
This completes the proof. 
\qed 

\begin{prop}\label{resolvent}
There exist positive constants $\nu_0$, $\gamma_0$, $\varepsilon_0$ and $a$ such that if $\frac{\nu^2}{\nu+\tilde{\nu}} \geq \nu_0$, $\frac{\gamma^2}{\nu+ \tilde{\nu}} \geq \gamma_0$ and $S\leq \varepsilon_{0}\frac{\nu^2}{\gamma^2\sqrt{\nu+\tilde{\nu}}}\sqrt{1-e^{-a\frac{\nu+\tilde{\nu}}{\gamma^2}}}$, then 
\begin{equation*}
\Sigma_0:=\{\lambda:\mathrm{Re}\lambda > -\frac{\beta_0}{2};~\lambda\neq2\pi ik, \forall k \in \mathbb{Z}\} \subset \rho(-B_0).
\end{equation*}
Furthermore, $u=\trans(\phi,w)=(\lambda+B_0)^{-1}F$ satisfies the estimates 
\begin{align*}
\| u\|_{X}&\leq \frac{C}{(2\mathrm{Re}\lambda+\beta_0)^\frac{1}{2}|1-e^{-\lambda}|}\int^1_0|\langle f(s)\rangle |\,ds \\
&\hspace{4mm}+\frac{C}{(1-e^{-(2\mathrm{Re}\lambda+\beta_0)})^\frac{1}{2}}\|F\|_{L^2(\mathbb{T}_1;L^2(\Omega_{per})\times H^{-1}(\Omega_{per}))}, \\
\|\nabla w\|_{L^2(\mathbb{T}_1;L^2(\Omega_{per}))}&\leq \frac{C}{(2\mathrm{Re}\lambda+\beta_0)^\frac{1}{2}|1-e^{-\lambda}|}\int^1_0|\langle f(s)\rangle |\,ds \\
&\hspace{4mm}+\frac{C}{(1-e^{-(2\mathrm{Re}\lambda+\beta_0)})^\frac{1}{2}}\|F\|_{L^2(\mathbb{T}_1;L^2(\Omega_{per})\times H^{-1}(\Omega_{per}))} 
\end{align*}
\end{prop}

We are now in a position to prove Theorem \ref{critical eigen}.

\vspace{2ex} 

\hspace{-6mm}{\bf Proof of Theorem \ref{critical eigen}.} 
We first observe that 
\begin{gather}\label{est b1}
\begin{split}
&\|B_j^{(1)}u\|_{L^2([0,1];L^2(\Omega_{per}) \times L^2(\Omega_{per}))} \\
&\hspace{8mm}\leq C\left(\|u\|_{L^2([0,1];L^2(\Omega_{per})\times L^2(\Omega_{per}))}+ \|\nabla w\|_{L^2([0,1];L^2(\Omega_{per}))}\right)  \end{split} 
\end{gather}
and 
\begin{gather}\label{est b2}
\|B_{j,k}^{(2)}u\|_{L^2([0,1];L^2(\Omega_{per}) \times L^2(\Omega_{per}))} \leq C\|w\|_{L^2([0,1];L^2(\Omega_{per}))}.
\end{gather}

Let $\Sigma$ be the set given in Theorem \ref{critical eigen}. 
We see from Proposition \ref{resolvent} that if $\lambda\in\Sigma$, then 
\begin{gather*}
\|(\lambda+B_0)^{-1}F\|_{L^2(\mathbb{T}_1;L^2(\Omega_{per})\times H^1(\Omega_{per}))}\leq C\| F\|_X. 
\end{gather*}
This, together with (\ref{est b1}) and (\ref{est b2}), implies that 
\begin{gather*}
\|M_{\eta'}(\lambda+B_0)^{-1}F\|_X\leq C(|\eta'|+|\eta'|^2)\| F\|_X. 
\end{gather*}
It then follows that there exists positive constant $r_0=r_0(\nu, \tilde{\nu}, \gamma, \nu_0, \gamma_0)$ such that if $|\eta'|\leq r_0$, then 
\begin{gather*}
\|M_{\eta'}(\lambda+B_0)^{-1}F\|_X\leq \frac{1}{2}\| F\|_X.
\end{gather*}
We thus find that if $|\eta'|\leq r_0$, then $\Sigma\subset \rho(-B_{\eta'})$ and for $\lambda\in\Sigma$ 
\begin{gather*}
(\lambda+B_{\eta'})^{-1}=(\lambda+B_0)^{-1}\sum^\infty_{N=0}(-1)^N[M_{\eta'}(\lambda+B_0)^{-1}]^N
\end{gather*}
and 
\begin{gather*}
\|(\lambda+B_{\eta'})^{-1}F\|_{L^2(\mathbb{T}_1;L^2(\Omega_{per})\times H^1(\Omega_{per}))}\leq 2C\| F\|_X. 
\end{gather*}
This proves the assertion (i).

As for the assertion (ii), it suffices to show that if $|\eta'|\leq r_0$, then 
\begin{gather*}
\sigma(-B_{\eta'})\cap \left\{\lambda\in\mathbb{C};\, \mathrm{Re}\lambda\geq-\frac{\beta_0}{4},\, |\lambda|\leq\frac{\pi}{4}\right\}=\{\lambda_{\eta',0}\}, \\
\lambda_{\eta',0}= -i\sum_{j=1}^{n-1}a_j\eta_j - \sum_{j,k=1}^{n-1}a_{jk}\eta_j \eta_k  + O(|\eta'|^3)  \hspace{4mm} (\eta' \rightarrow 0)
\end{gather*}
with some constants $a_j,a_{jk} \in \mathbb{R}$ and $\lambda_{\eta',0}$ satisfies  
\begin{gather*}
\mathrm{Re}\lambda_{\eta',0}\leq -\frac{\kappa_0\gamma^2}{2\nu}|\eta'|^2. 
\end{gather*}

In view of Proposition \ref{eigenspace}, Proposition \ref{resolvent}, (\ref{est b1}) and (\ref{est b2}), 
we can apply the analytic perturbation theory (\cite{K}) to see that the set 
\begin{equation*}
\sigma(-B_{\eta'})\cap \left\{\lambda\in\mathbb{C};\, \mathrm{Re}\lambda\geq-\frac{\beta_0}{4},\, |\lambda|\leq\frac{\pi}{4}\right\}
\end{equation*}
 consists of a simple eigenvalue, say $\lambda_{\eta',0}$, for sufficiently small $\eta'$, and that $\lambda_{\eta',0}$ is expanded as 
\begin{equation*}
\lambda_{\eta',0}= \lambda_0 + \sum_{j=1}^{n-1}\eta_j \lambda_j^{(1)}+ \sum_{j,k=1}^{n-1}\eta_j \eta_k \lambda_{jk}^{(2)} + O(|\eta'|^3)  \hspace{4mm} (\eta' \rightarrow 0),
\end{equation*}
where 
\begin{gather*}
\lambda_0=0, \\
\lambda_j^{(1)}=-\inn{B_j^{(1)}u^{(0)}}{u^{(0)*}}, \\
\lambda_{jk}^{(2)}=-\frac{1}{2}\inn{(B_{j,k}^{(2)}+B_{k,j}^{(2)})u^{(0)}}{u^{(0)*}} +\frac{1}{2}\inn{(B_j^{(1)}SB_k^{(1)}+B_k^{(1)}SB_j^{(1)})u^{(0)}}{u^{(0)*}}.
\end{gather*}
Here $S=[(I-\Pi^{(0)})B_0(I-\Pi^{(0)})]^{-1}$. 
By definition of $B_j^{(1)}$, $u^{(0)}$ and $u^{(0)*}$, we have 
\begin{equation*}
\lambda_j^{(1)}=-i\mean{v_p^j\phi^{(0)}+\gamma^2 \rho_p w^{(0),j}}.
\end{equation*}
As for $\lambda_{jk}^{(2)}$, since $\frac{1}{2}\inn{(B_{j,k}^{(2)}+B_{k,j}^{(2)})u^{(0)}}{u^{(0)*}}=0$, we obtain
\begin{equation*}
\lambda_{jk}^{(2)}=\frac{1}{2}\inn{(B_j^{(1)}SB_k^{(1)}+B_k^{(1)}SB_j^{(1)})u^{(0)}}{u^{(0)*}}.
\end{equation*}
We set $u_1^{(k)}=\trans(i\phi_1^{(k)},\frac{i}{\nu}w_1^{(k)})=SB_k^{(1)}u^{(0)}$. Then $u_1^{(k)}$ is a solution of 
\begin{equation} \label{I-Pi}
\begin{cases}
B_0u_1=(I-\Pi^{(0)})B_k^{(1)}u^{(0)},\\
\mean{\phi_1}=0,
\end{cases}
\end{equation}
where 
\begin{align*}
&(I-\Pi^{(0)})B_k^{(1)}u^{(0)} \\
&=i\begin{pmatrix}
v_p^k\phi^{(0)} + \gamma^2 \rho_p w^{(0)}\cdot \mathbf{e}_k \\
(\frac{P'(\rho_p)}{\gamma^2\rho_p}\phi^{(0)})\mathbf{e}_k-\frac{1}{\rho_p}(2\nu \partial_{x_k} w^{(0)} + \mathbf{e}_k\tilde{\nu} \mathrm{div}w^{(0)}- \tilde{\nu} \nabla w^{(0),k}) 
\end{pmatrix} \\
&\hspace{4mm}-i \ll v_p^k\phi^{(0)} + \gamma^2 \rho_p w^{(0)}\cdot \mathbf{e}_k \gg u^{(0)}. 
\end{align*}
In fact, there exists a solution $u_1^{(k)}=\trans(\phi_1^{(k)}, w_1^{(k)})$ to (\ref{I-Pi}), and $\lambda^{(2)}_{jk}$ is written as 
\begin{align*}
\lambda^{(2)}_{jk}
&=\frac{1}{2}\inn{B_j^{(1)}u_1^{(k)}+B_k^{(1)}u_1^{(j)}}{u^{(0)*}} \\
&=\frac{1}{2}\mean{v_p^j\phi_1^{(k)} + \frac{\gamma^2 }{\nu}\rho_p w_1^{(k)}\cdot\mathbf{e}_j}
+\frac{1}{2}\mean{v_p^k\phi_1^{(j)} + \frac{\gamma^2 }{\nu}\rho_p w_1^{(j)}\cdot\mathbf{e}_k}.
\end{align*}
We thus estimate $u_1^{(k)}$ to prove the estimate (\ref{eigen}). 
\begin{lem}\label{lem4.9}
Assume that $\frac{\nu^2}{\nu+\tilde{\nu}} \geq \nu_0$, $\frac{\gamma^2}{\nu+ \tilde{\nu}} \geq \gamma_0$, $S\leq \varepsilon_{0}\frac{\nu^2}{\gamma^2\sqrt{\nu+\tilde{\nu}}}\sqrt{1-e^{-a\frac{\nu+\tilde{\nu}}{\gamma^2}}}$ and $|\eta'| \leq r_0 $. 
Then the following estimate holds$:$
\begin{gather} \label{eigen1}
\int_{0}^{1} \left( \delta \|\phi_1^{(k)}(s) \|_{L^2(\Omega_{per})}^2 + \frac{1}{\nu}\|\nabla w_1^{(k)}(s)\|_{L^2(\Omega_{per})}^2\right)\,ds 
\leq \frac{C}{\nu}. 
\end{gather}
\end{lem}

\hspace{-6mm}{\bf Proof.}
As in Proposition \ref{L2 time}, we have 
\begin{align*} 
&\int^1_0\left( \delta \|\phi_1^{(k)}(s) \|_{L^2(\Omega_{per})}^2 + \frac{1}{\nu}\|\nabla w_1^{(k)}(s)\|_{L^2(\Omega_{per})}^2\right) ds \\
&\leq C\left(\frac{1}{\delta\gamma^4}+\frac{\delta^2}{\nu} \right)\int^1_0\left\|\left((I-\Pi^{(0)})B_k^{(1)}u^{(0)}\right)_1\right\|_{L^2(\Omega_{per})}^2~ds \\
&\hspace{4mm}+\frac{C}{\nu}\int^1_0\left\|\left((I-\Pi^{(0)})B_k^{(1)}u^{(0)} \right)_2\right\|^2_{H^{-1}(\Omega_{per})}\,ds, 
\end{align*} 
where $(u)_1=\phi$ and $(u)_2=w$ for $u=\trans(\phi, w)$. 
Since $\frac{1}{\delta\gamma^4}+\frac{\delta^2}{\nu}\leq \frac{1}{\gamma^2}\left(\frac{1}{\nu}+\frac{1}{\gamma}+\frac{\nu+\tilde{\nu}}{\gamma^2}\right)$, it follows from Propositions \ref{timeper} and \ref{exist} that 
\begin{align}
\left(\frac{1}{\delta\gamma^4}+\frac{\delta^2}{\nu}\right)\int^1_0\left\|\left((I-\Pi^{(0)})B_k^{(1)}u^{(0)}\right)_1\right\|_{L^2(\Omega_{per})}^2\,ds &\leq \frac{C}{\nu}\left(\frac{1}{\nu^2}+\frac{1}{\nu\gamma}+\frac{1}{\gamma^2} \right),  \label{B_11}\\ 
\frac{1}{\nu}\int^1_0\left\|\left((I-\Pi^{(0)})B_k^{(1)}u^{(0)} \right)_2\right\|^2_{H^{-1}(\Omega_{per})}\,ds &\leq \frac{C}{\nu}, \notag
\end{align}
and we have the desired estimate. 
\qed 

\vspace{2ex}

To derive the estimate (\ref{eigen}), we next introduce $\tilde{u}_1^{(k)}=\trans(i\tilde{\phi}_1^{(k)},\frac{i}{\nu}\tilde{w}_1^{(k)})$, which is a unique stationary solution of the Stokes system
  \begin{equation} \label{Stokes}
 A\tilde{u}_1^{(k)}=F^{(k)}, ~ \langle\tilde{\phi}_1^{(k)}\rangle=0, 
  \end{equation}
  where 
\begin{equation*}
F^{(k)}=
\begin{pmatrix}
0\\ \mathbf{e}_k
\end{pmatrix}.
\end{equation*}
We use the following lemma (\cite[Theorem 4.7]{Makio}).
\begin{lem}[\cite{Makio}]\label{lem Makio}
Let $\tilde{\kappa}(\eta')$ be defined by 
\begin{equation*}
\tilde{\kappa}(\eta')=\sum^{n-1}_{j,k=1}\tilde{a}_{jk}\eta_j\eta_k 
\end{equation*}
for $\eta'\in\mathbb{R}^{n-1}$, where 
\begin{equation*}
\tilde{a}_{jk}=\frac{\gamma^2}{\nu}(\nabla \tilde{w}_1^{(j)},\nabla \tilde{w}_1^{(k)})=\frac{\gamma^2}{\nu}\langle \mathbf{e}_j\cdot\tilde{w}_1^{(k)}\rangle.
\end{equation*}
Then there exists a constant $\kappa_0>0$ independent of $\nu$, $\tilde{\nu}$ and $\gamma$ such that 
\begin{equation*}
\tilde{\kappa}(\eta')\geq \frac{\kappa_0\gamma^2}{\nu}|\eta'|^2 
\end{equation*}
for all $\eta'\in\mathbb{R}^{n-1}$.
\end{lem}

By using Lemma \ref{lem Makio}, we have the following estimate. 

\begin{lem}\label{lem4.10}
Assume that $\frac{\nu^2}{\nu+\tilde{\nu}} \geq \nu_0$, $\frac{\gamma^2}{\nu+ \tilde{\nu}} \geq \gamma_0$, $S\leq \varepsilon_{0}\frac{\nu^2}{\gamma^2\sqrt{\nu+\tilde{\nu}}}\sqrt{1-e^{-a\frac{\nu+\tilde{\nu}}{\gamma^2}}}$ and $|\eta'| \leq r_0 $, then the following estimate holds$:$
\begin{gather} 
\int^1_0 \|\nabla(w_1^{(k)}-\tilde{w}_1^{(k)})\|_{L^2(\Omega_{per})}^2\,ds
\leq C\left(\frac{1}{\nu^2}+\frac{1}{\nu\gamma}+\frac{1}{\gamma^2} \right). \label{eigen2}
\end{gather}
   \end{lem}
\hspace{-6mm}{\bf Proof.}
We consider 
\begin{equation*}
(\partial_t+L(t))(u_1^{(k)}-\tilde{u}_1^{(k)})=(I-\Pi^{(0)})B_k^{(1)}u^{(0)}-F^{(k)}-M\tilde{u}_1^{(k)}.
\end{equation*}
It follows from the estimate for the Stokes problem (see, e.g., \cite{G}) that $\|\partial_x\tilde{\phi}_1^{(k)}\|^2_{L^2(\Omega_{per})}+\|\partial_x^2\tilde{w}_1^{(k)}\|^2_{L^2(\Omega_{per})}\leq C$, and we have 
\begin{equation} \label{est M}
\| M\tilde{u}_1^{(k)}\|^2_{L^2(\mathbb{T}_1; L^2(\Omega_{per})\times H^{-1}(\Omega_{per}))}\leq C\frac{\nu+\tilde{\nu}}{\gamma^4}.
\end{equation}
By using (\ref{B_11}), (\ref{est M}) and 
\begin{equation}
\int^1_0\left\|\left((I-\Pi^{(0)})B_k^{(1)}u^{(0)}-F^{(k)} \right)_2\right\|^2_{H^{-1}(\Omega_{per})}\,ds \leq C\left(\frac{1}{\gamma^2}+\frac{1}{\nu^4} \right), 
\end{equation}
we obtain 
\begin{equation*}
\frac{1}{\nu}\int^1_0\|\nabla(w_1^{(k)}-\tilde{w}_1^{(k)})\|_{L^2(\Omega_{per})}^2\,ds\leq \frac{C}{\nu}\left(\frac{1}{\nu^2}+\frac{1}{\nu\gamma}+\frac{1}{\gamma^2} \right). 
\end{equation*}
This completes the proof. \qed

\vspace{2ex}

 \hspace{-6mm}{\bf Proof of (\ref{eigen}).}
 By  Lemmas \ref{lem4.9}, \ref{lem Makio} and \ref{lem4.10}, if $\frac{\nu^2}{\nu+\tilde{\nu}} \geq \nu_0$, $\frac{\gamma^2}{\nu + \tilde{\nu}}\geq \gamma_0$ and $S\leq \varepsilon_{0}\frac{\nu^2}{\gamma^2\sqrt{\nu+\tilde{\nu}}}\sqrt{1-e^{-a\frac{\nu+\tilde{\nu}}{\gamma^2}}}$, we have 
  \begin{align*}
  \sum_{j,k=1}^{n-1}\eta_j \eta_k \lambda_{jk}^{(2)}&=-\sum_{j,k=1}^{n-1}\eta_j \eta_k \mean{v_p^j\phi_1^{(k)}(t) + \frac{\gamma^2 }{\nu}\rho_p w_1^{(k)}(t)\cdot\mathbf{e}_j}\\
  &=-\sum_{j,k=1}^{n-1}\eta_j \eta_k \mean{v_p^j\phi_1^{(k)}(t) + \frac{\gamma^2 }{\nu}\phi_p w_1^{(k)}(t)\cdot \mathbf{e}_j}\\
   &\hspace{4mm}-\sum_{j,k=1}^{n-1}\eta_j \eta_k (\mean{(w_1^{(k)}(t)-\tilde{w}_1^{(k)}(t))\cdot \mathbf{e}_j }+\frac{\gamma^2}{\nu}\mean{\tilde{w}_1^{(k)}(t)\cdot \mathbf{e}_j })\\
   &\leq -\frac{\gamma^2}{\nu}\left(\kappa_0-\frac{\nu}{\gamma^2}\right) |\eta'|^2 \\
&\leq - \frac{\kappa_0\gamma^2}{2\nu}|\eta'|^2.
  \end{align*}
  This complets the proof. \qed 

\vspace{2ex}

We next prove Theorem \ref{u_eta}. 
To do so, we establish the estimate for $(\lambda + B_0)^{-1}F$ in a higher order Sobolev space.
\begin{lem} \label{higher}
There exist positive constants $\nu_0$, $\gamma_0$, $\varepsilon_0$ and $a$ such that if  $\frac{\nu^2}{\nu+\tilde{\nu}} \geq \nu_0$, $\frac{\gamma^2}{\nu + \tilde{\nu}}\geq \gamma_0,$ and $S\leq \varepsilon_{0}\frac{\nu^2}{\gamma^2\sqrt{\nu+\tilde{\nu}}}\sqrt{1-e^{-a\frac{\nu+\tilde{\nu}}{\gamma^2}}}$ for $\lambda$ satisfying uniformly $\frac{\pi}{4} \leq |\lambda| \leq \frac{\pi}{2}$, then $u=\trans(\phi, w)=(\lambda+B_0)^{-1}F$ with $F=\trans(f, g)$ satisfies
 \begin{equation}
\nom{u}{2}{2}
+\int^t_0\left(\nom{\phi(s)}{2}{2}+\nom{w(s)}{2}{3}\right)\,ds 
\leq C \int_0^1 \left(\nom{f(s)}{2}{2}+\nom{g(s)}{2}{1}\right)\, ds.
 \end{equation}
 \end{lem}
 \hspace{-6mm}{\bf Proof.} 
We consider
\begin{equation}
(\lambda + B_0)u= F,\hspace{4mm}u\in D(B_0).  \label{pi}
\end{equation} 
We set $\Pi_1= I - \Pi^{(0)}$. 
Applying $\Pi^{(0)}$ and $\Pi_1$, (\ref{pi}) is decomposed into 
  \begin{gather}\label{Pi1}
  \lambda \Pi^{(0)}u = \Pi^{(0)}F, \\
    \lambda u_1 + B_0u_1  = G_1, \label{Pi2}
    \end{gather}
where $u_1=\Pi_1 u$ and $G_1=\trans(g^0_1, \tilde{g}_1):=\Pi_1F$. 
If $\lambda \neq 0$, then we rewrite (\ref{Pi1}) as 
  \begin{equation} \label{Pi3}
  \Pi^{(0)}u=\frac{1}{\lambda}\Pi^{(0)}F.
  \end{equation}
We next consider (\ref{Pi2}). 
We see from the proofs of Propositions \ref{eigenspace} and \ref{prop C} that there exists a unique solution $u_1\in D(B_0)\cap X_1$ to (\ref{Pi2}) if $\lambda\in \Sigma_0$. 
Furthermore, it follows from the proofs of Lemma \ref{H2 global} and Proposition \ref{prop C} that 
\begin{align*} 
&\nom{u_1}{2}{2} + \int^t_0e^{-(2\mathrm{Re}\lambda+a\frac{\nu+\tilde{\nu}}{\gamma^2})(t-s)}\left(\nom{\phi_1(s)}{2}{2}+\nom{w_1(s)}{2}{3}\right)\,ds \\
&\leq \frac{C}{1-e^{-(2\mathrm{Re}\lambda+a\frac{\nu+\tilde{\nu}}{\gamma^2})}}\int_0^1\left(\nom{g^0_1(s)}{2}{2}+\nom{\tilde{g}_1(s)}{2}{1}\right)\,ds. 
\end{align*}
Let $\frac{\pi}{4} \leq |\lambda| \leq \frac{\pi}{2}$. 
We set $$u=\frac{1}{\lambda}\Pi^{(0)}F+u_1.$$ 
Then $u$ is a solution of $(\lambda + B_0)u=G$, where $G=\trans(g^0, \tilde{g})=\Pi^{(0)}F+G_1$. Then it holds that 
\begin{equation} \label{u_1_4}
\begin{split}
&\nom{u}{2}{2}+\int^t_0\left(\nom{\phi(s)}{2}{2}+\nom{w(s)}{2}{3}\right)\,ds \\ 
&\leq C\left(\frac{1}{|\lambda|}+\frac{1}{1-e^{-(2\mathrm{Re}\lambda+a\frac{\nu+\tilde{\nu}}{\gamma^2})}} \right)
\int^1_0\left(\nom{g^0(s)}{2}{2}+\nom{\tilde{g}(s)}{2}{1}\right)\,ds.
\end{split}
\end{equation}
     This implies that $(\lambda + B_0)^{-1}$ is bounded and $u=(\lambda + B_0)^{-1}G$ satisfies \eqref{u_1_4}. \\
Since 
\begin{equation*}
\int^1_0\left(\nom{g^0(s)}{2}{2}+\nom{\tilde{g}(s)}{2}{1}\right)\,ds \leq C\int^1_0\left(\nom{f(s)}{2}{2}+\nom{g(s)}{2}{1}\right)\, ds, 
\end{equation*}  
there exist positive constants $\nu_0$, $\gamma_0$ and $\varepsilon_0$ such that if $\frac{\nu^2}{\nu+\tilde{\nu}} \geq \nu_0$, $\frac{\gamma^2}{\nu+\tilde{\nu}}\geq \gamma_0$ and $S\leq \varepsilon_{0}\frac{\nu^2}{\gamma^2\sqrt{\nu+\tilde{\nu}}}\sqrt{1-e^{-a\frac{\nu+\tilde{\nu}}{\gamma^2}}}$, then $\lambda + B_0$ has a bounded inverse for $\lambda$ satisfying $\frac{\pi}{4} \leq |\lambda| \leq \frac{\pi}{2}$ and $u=\trans(\phi, w)=(\lambda+B_0)^{-1}F$ satisfies 
\begin{align*}
       \nom{u}{2}{2} + \int^t_0\left(\nom{\phi(s)}{2}{2}+\nom{w(s)}{2}{3}\right)\,ds 
\leq C \int_0^1 \left(\nom{f(s)}{2}{2}+\nom{g(s)}{2}{1}\right)\, ds.
        \end{align*}
This completes the proof.   \qed 

\vspace{2ex}

We are now in a position to prove Theorem \ref{u_eta}. 

\vspace{2ex}

\hspace{-6mm}{\bf Proof of Theorem \ref{u_eta}.}
If $\lambda \in \Sigma$, then for $|\eta'| \leq r_0$, $(\lambda + B_{\eta'})^{-1}$ is given by the Neumann series expansion 
\begin{equation*}
(\lambda + B_{\eta'})^{-1}=(\lambda + B_{0})^{-1} \sum_{N=0}^{\infty}(-1)^N[M_{\eta'}(\lambda + B_{\eta'})^{-1}]^N.
\end{equation*}
It then follows that 
\begin{gather*}
u_{\eta'}=u^{(0)}+i\eta' \cdot u^{(1)} + u^{(2)}, \\
u_{\eta'}^*=u^{(0)*}+i\eta' \cdot u^{(1)*} + u^{(2)*}, \\
\langle \tilde{u}_{\eta'}, \tilde{u}_{\eta'}^* \rangle_t = 1 + O(\eta') \geq \frac{1}{2} 
\end{gather*}
for $|\eta'| \leq r_0$.
Here, 
\begin{equation*}
u^{(1)}=- \frac{1}{2\pi i}\int_{|\lambda|=\frac{\pi}{4}}(\lambda + B_{0})^{-1}B_j^{(1)}(\lambda + B_{0})^{-1}u^{(0)}\,d\lambda,
\end{equation*}
\begin{equation*}
u^{(2)}(\eta')=\frac{1}{2\pi i}\int_{|\lambda|=\frac{\pi}{4}}R^{(2)}(\lambda, \eta)u^{(0)}\,d\lambda,
\end{equation*}
with
\begin{multline*}
R^{(2)}(\lambda,\eta)=-(\lambda + B_{0})^{-1} B_{j,k}^{(2)}(\lambda + B_{0})^{-1}\\
+(\lambda + B_{0})^{-1}\sum_{N=2}^{\infty}(-1)^N \eta^{N-2}[(B_j^{(1)} + \eta B_{j,k}^{(2)})(\lambda + B_{0})^{-1}]^{N}.
\end{multline*}
By Proposition \ref{higher}
and the definition of $B_j^{(1)} ,B_{j,k}^{(2)}$, we obtain 
\begin{align*} 
\|u^{(1)}\|_{C([0,1];H^2(\Omega_{per}))} 
\leq C
\quad \text{and} \quad 
  \|u^{(2)}(\eta)\|_{C([0,1];H^2(\Omega_{per}))} \leq  C.
  \end{align*} 
Consequently, we have  
\begin{gather*}
\|u_{\eta'}(t)\|_{H^2(\Omega_{per})} \leq C, \\
\no{u_{\eta'}(t)-u^{(0)}(t)}{}{H^2(\Omega_{per})}\leq C|\eta'|. 
\end{gather*}
 Similarly, we obtain 
\begin{gather*}
\|\tilde{u}_{\eta'}^*(t)\|_{H^2(\Omega_{per})} \leq C.
\end{gather*}
This completes the proof. \qed

\section{Proof of Proposition \ref{timeper}}\label{sec5}

In this section we give a proof of Proposition \ref{timeper}.
We set 
\begin{equation*}
\rho=1+\frac{1}{\gamma^2}\phi, \quad v=w. 
\end{equation*} 
The system (\ref{dim1})-(\ref{dim2}) is then written as 
\begin{gather} \label{prob3}
\partial_t \phi + \gamma^2 \mathrm{div} w = f(\phi,w),  \\
\label{prob2}
\partial_t w -\nu \Delta w- \tilde{\nu}\nabla \mathrm{div} w + \nabla \phi = g(\phi,w,G), 
\end{gather}
where 
\begin{gather*}
f(\phi,w)=-\mathrm{div}(\phi w), \\
\begin{split}
g(\phi,w,G)=&SG-w\cdot\nabla w-\frac{\phi}{\gamma^2+\phi}\{\nu\Delta w+\tilde{\nu}\nabla\mathrm{div}w\} \\
&+\frac{\phi}{\gamma^2+\phi}\nabla\phi-\frac{1}{\gamma^2+\phi}\nabla\left(p^{(1)}\left(\frac{\phi}{\gamma^2}\right)\phi^2\right).
\end{split}
\end{gather*}
We consider (\ref{prob3})-(\ref{prob2}) on $\Omega_{per}$ 
under the conditions 
\begin{gather}
w|_{x_n=0,1}=0, \label{BC} \\
\langle \phi \rangle  =0. \label{meanvalue}
\end{gather}

Under some smallness assumption on the size of $S$, 
we have the following result on the existence of a time periodic solution of (\ref{prob3})-(\ref{meanvalue}).
\begin{prop}\label{prop1}
Let $G\in \cap_{j=0}^2H^j(\mathbb{T}_1;H^{3-2j}(\Omega_{per}))$ with $[G]_{3,1,\Omega_{per}}=1$.
There exist positive constants $\nu_0$, $\gamma_0$, $\varepsilon_0$ and $a$ such that 
if $\frac{\nu^2}{\nu+\tilde{\nu}}\geq \nu_0$, $\frac{\gamma^2}{\nu+\tilde{\nu}}\geq \gamma_0$ and $S\leq \varepsilon_{0}\frac{\nu^2}{\gamma^2\sqrt{\nu+\tilde{\nu}}}\sqrt{1-e^{-a\frac{\nu+\tilde{\nu}}{\gamma^2}}}$, then there exists a time periodic solution $u=(\phi,w)\in \cap_{j=0}^2C^j(\mathbb{T}_1;H^{4-2j}(\Omega_{per})\times H^{4-2j}(\Omega_{per}))\cap H^j(\mathbb{T}_1;H^{4-2j}(\Omega_{per})\times H^{5-2j}(\Omega_{per}))$ to problem $(\ref{prob3})$-$(\ref{meanvalue})$, and $u$ satisfies
\begin{gather*}
\sup_{t\in\mathbb{T}_1}\left\{\frac{1}{\gamma^2}[\![\phi(t)]\!]_4^2+[\![w(t)]\!]_4^2\right\}
\leq C\frac{\nu^2}{\gamma^4}, \\
\int_{\mathbb{T}_1}\left\{\frac{\nu^2}{\nu+\tilde{\nu}}[\![ w(t) ]\!]^2_5+\frac{1}{\nu+\tilde{\nu}}[\![ \nabla\phi(t) ]\!]^2_3 + \frac{\nu+\tilde{\nu}}{\gamma^4}\|\partial_t^2\phi(t)\|^2_{L^2(\Omega_{per})}\right\}\,dt\leq C\frac{\nu^2}{\gamma^4},
\end{gather*}
where $C$ is independent of $\nu$, $\tilde{\nu}$, $\gamma$ and $S$.
\end{prop}

Let $u= \trans(\phi, w)$ be the time periodic solution obtained in Proposition \ref{prop1}. 
Setting $\phi_{p}=\frac{1}{\gamma^{2}}\phi$, $v_{p} = w$ and $u_{p}=\trans(1+\phi_{p}, v_{p})$, we obtain Proposition \ref{timeper}.

\begin{rem}
If $\frac{\gamma^2}{\nu+\tilde{\nu}}\gg1$, then $1-e^{-a\frac{\nu+\tilde{\nu}}{\gamma^2}}\sim a\frac{\nu+\tilde{\nu}}{\gamma^2}$, and so the condition of $S$ in Proposition $\ref{prop1}$ implies  
\begin{equation*}
S\leq \varepsilon_{0}\frac{\nu^2}{\gamma^2\sqrt{\nu+\tilde{\nu}}}\sqrt{1-e^{-a\frac{\nu+\tilde{\nu}}{\gamma^2}}}
\sim \varepsilon_0\sqrt{a}\frac{\nu^2}{\gamma^3}. 
\end{equation*}

\end{rem}

Proposition \ref{prop1} follows from Valli's argument (see \cite{Valli}). 
Since we look for a solution in a higher order Sobolev space than that considered in \cite{Valli} and we need to take care of the dependence of the estimates on the parameters, we here give a proof of Proposition \ref{prop1}

To prove Proposition \ref{prop1} we first consider the initial boundary problem for (\ref{prob3})-(\ref{prob2}) under the boundary condition (\ref{BC}) and the initial condition 
\begin{equation}
u|_{t=0}=u_0=\trans(\phi_0,w_0). \label{IC}
\end{equation}

In what follows we assume that 
\begin{equation}\label{cc0} 
\begin{split}
&u_0\in H^4_*(\Omega_{per})\times (H^4(\Omega_{per})\cap H^1_0(\Omega_{per})), \\ 
&G \in \cap_{j=0}^2H^j(\mathbb{T}_1;H^{3-2j}(\Omega_{per})), \quad [G]_{3,1,\Omega_{per}} = 1.
\end{split}
\end{equation}
We will also impose the following compatibility conditions on $u_0=\trans(\phi_0,w_0)$ and $G$:
\begin{gather}
w_0\in H^1_0(\Omega_{per}), \label{cc1}\\
\nu\Delta w_0+\tilde{\nu}\nabla\mathrm{div}w_0-\nabla\phi_0+g(\phi_0,w_0,G(0))\in H^1_0(\Omega_{per}). \label{cc2}
\end{gather}
We note that if $u_0\in H^4_*(\Omega_{per})\times\left( H^4(\Omega_{per})\cap H^1_0(\Omega_{per})\right)$ and $G\in \cap_{j=0}^2H^j(\mathbb{T}_1;H^{3-2j}(\Omega_{per}))$, then $\nu\Delta w_0+\tilde{\nu}\nabla\mathrm{div}w_0-\nabla\phi_0+g(\phi_0,w_0,G(0))\in H^1_0(\Omega_{per})$.

\subsection{Global solution}\label{sec5.1}
In this subsection, we prove the global existence of solution to (\ref{prob3})-(\ref{IC}) with energy estimate that is stated in the following proposition. 

\begin{prop} \label{est}
Let $u_{0}$ and $G$ satisfy \eqref{cc0} and the compatibility conditions $(\ref{cc1})$ and $(\ref{cc2})$.
There exist positive constants $\nu_0$, $\gamma_0$, $a$, $C_0$, $C_1$ and $C_E$ such that 
if $\frac{\nu^2}{\nu+\tilde{\nu}}\geq \nu_0$, $\frac{\gamma^2}{\nu+\tilde{\nu}}\geq \gamma_0$, $S\leq\frac{1}{10 \sqrt{C_0(1 + C_E)}}\frac{\nu}{\sqrt{\nu+\tilde{\nu}}}\sqrt{1-e^{-a\frac{\nu+\tilde{\nu}}{\gamma^2}}}\min\left\{1, \frac{1}{C_1}\frac{\nu}{\gamma} \right\}$ and $E_4(u_0)\leq \frac{4C_0(1+C_E)}{\nu}S^2\frac{1}{1-e^{-a\frac{\nu+\tilde{\nu}}{\gamma^2}}}$, then there exists a unique global solution $u$ to $(\ref{prob3})$-$(\ref{IC})$ in $\cap_{j=0}^2C^j([0,\infty);H^{4-2j}(\Omega_{per})\times H^{4-2j}(\Omega_{per}))\cap H^j([0,\infty);H^{4-2j}(\Omega_{per})\times H^{5-2j}(\Omega_{per}))$ which satisfies  
\begin{equation}\label{H4}
\begin{split} 
& E_4(u(t))+\frac{1}{4}\int^t_0e^{-a\frac{\nu+\tilde{\nu}}{\gamma^2}(t-s)}D_4(u(s))\,ds \\[1ex]
& \ \ \leq (1+C_E) \left\{ e^{-a\frac{\nu+\tilde{\nu}}{\gamma^2}t}E_4(u_0)+2\frac{C_0}{\nu}S^2\frac{1}{1-e^{-a\frac{\nu+\tilde{\nu}}{\gamma^2}}} \right\},
\end{split}
\end{equation}
where $E_m(u)$ and $D_m(u)$ $(m=2,4)$ are the quantities satisfying the following inequalities with some positive constant $C:$  

\begin{gather*}
C^{-1}\left(\frac{1}{\gamma^2}[\![ \phi ]\!]^2_m+[\![w]\!]_m^2\right) \leq E_m(u) \leq C\left(\frac{1}{\gamma^2}[\![ \phi ]\!]^2_m+[\![w]\!]_m^2\right) \\
\begin{split}
&C^{-1}\left(\frac{\nu^2}{\nu+\tilde{\nu}}[\![ w ]\!]^2_{m+1}+\frac{1}{\nu+\tilde{\nu}}[\![ \nabla\phi ]\!]^2_{m-1} + \frac{\nu+\tilde{\nu}}{\gamma^4}\|\partial_t^{\frac{m}{2}}\phi\|^2_{L^2(\Omega_{per})}\right) \leq D_m(u) \\
&\leq C\left(\frac{\nu^2}{\nu+\tilde{\nu}}[\![ w ]\!]^2_{m+1}+\frac{1}{\nu+\tilde{\nu}}[\![ \nabla\phi ]\!]^2_{m-1} + \frac{\nu+\tilde{\nu}}{\gamma^4}\|\partial_t^{\frac{m}{2}}\phi\|^2_{L^2(\Omega_{per})}\right). 
\end{split}
\end{gather*}
\end{prop}

As in \cite{E}, we can prove Proposition \ref{est} by combining local existence and the a priori estimates. 
The local existence is proved by applying the local solvability result in \cite{Kagei-Kawashima2,Valli}.
In fact, we can show that the following assertion.

\begin{prop}\label{local}
Let $u_0$ and $G$ satisfy \eqref{cc0} and the compatibility conditions $(\ref{cc1})$ and $(\ref{cc2})$.
Then there exists a positive number $T$ depending only on $\| u_0\|_{H^4(\Omega_{per})}$, $\nu$, $\tilde{\nu}$, $\gamma$ and $S$ such that the problem $(\ref{prob3})$-$(\ref{IC})$ has a unique solution $u(t)=\trans(\phi(t),w(t))\in \cap_{j=0}^2C^j([0,T];H^{4-2j}(\Omega_{per})\times H^{4-2j}(\Omega_{per}))\cap H^j([0,T];H^{4-2j}(\Omega_{per})\times H^{5-2j}(\Omega_{per}))$ satisfying 
\begin{align*}
[\![u(t)]\!]_4^2+\int^T_0[\![\phi(s)]\!]_4^2+[\![w(s)]\!]_5^2\,ds \leq C\| u_0\|_{H^4(\Omega_{per})}^2.
\end{align*}
\end{prop}
 
The global existence is proved by combining Proposition \ref{local} and the following a priori estimates. 

\begin{prop}\label{apriori}
Let $T$ be a positive number and assume that $u$ is a solution of $(\ref{prob3})$-$(\ref{IC})$ in $\cap_{j=0}^2C^j([0,T];H^{4-2j}(\Omega_{per})\times H^{4-2j}(\Omega_{per}))\cap H^j([0,T]; H^{4-2j}(\Omega_{per})\times H^{5-2j}(\Omega_{per}))$. 
There exist positive constants $\nu_0$, $\gamma_0$, $a$, $C_0$, $C_1$ and $C_E$ independent of $T$, $\nu$, $\tilde{\nu}$, $\gamma$ and $S$ such that the following assertion holds. 
If $\frac{\nu^2}{\nu+\tilde{\nu}}\geq \nu_0$, $\frac{\gamma^2}{\nu+\tilde{\nu}}\geq \gamma_0$, $S\leq\frac{1}{10 \sqrt{C_0(1 + C_E)}}\frac{\nu}{\sqrt{\nu+\tilde{\nu}}}\sqrt{1-e^{-a\frac{\nu+\tilde{\nu}}{\gamma^2}}}\min\left\{1, \frac{1}{C_1}\frac{\nu}{\gamma} \right\}$ and $E_4(u(t))\leq \frac{6C_0(1+C_E)}{\nu}S^2\frac{1}{1-e^{-a\frac{\nu+\tilde{\nu}}{\gamma^2}}}$ for $t\in[0,T]$, then 
\begin{equation}\label{aprioriest}
\begin{split} 
& E_4(u(t))+ \frac{1}{4}\int_0^te^{-a\frac{\nu+\tilde{\nu}}{\gamma^2}(t-s)}D_4(u(s))\,ds  \\[1ex]
& \ \ \le (1+C_E) \left\{e^{-a\frac{\nu+\tilde{\nu}}{\gamma^2}t}E_4(u_0) 
+2\frac{C_0}{\nu}S^2\frac{1}{1-e^{-a\frac{\nu+\tilde{\nu}}{\gamma^2}}} \right\}. 
\end{split}
\end{equation}
for $t\in [0,T]$.
\end{prop}

Proposition \ref{apriori} is proved by the energy estimate and nonlinear estimates. 
In what follows, we denote 
\begin{align*}
&T_{j,l}= \partial_t^j \partial_{x'}^l.
\end{align*}
We also denote $\dot{\phi}:= \partial_t \phi + w \cdot \nabla \phi$. 
It then follows that 
 \begin{equation} \label{phidot}
 \dot{\phi} =-\gamma^2 \mathrm{div}w +\tilde{f},~~\tilde{f}=-\phi ~\mathrm{div}w.
 \end{equation}

We have the following basic estimates in a similar manner to \cite[Section 5]{E-K}. 

\begin{prop}\label{basicest}
The following estimates hold true$:$ 
\begin{equation} 
\begin{split}\label{1}
&\frac{1}{2}\frac{d}{dt}\no{T_{j,l}u}{2}{L^2(\Omega_{per}),\gamma}
+\nu\|T_{j,l}\nabla w\|_{L^2(\Omega_{per})}^2+\tilde{\nu}\|T_{j,l}\mathrm{div}w\|_{L^2(\Omega_{per})}^2+\frac{\nu + \tilde{\nu}}{\gamma^4}\|T_{j,l}\dot{\phi}\|_{L^2(\Omega_{per})}^2 \\
& \quad \leq (T_{j,l}F,T_{j,l}u)_\gamma  + C\frac{\nu + \tilde{\nu}}{\gamma^4}\|T_{j,l}\tilde{f}\|_{L^2(\Omega_{per})}^2,  
\end{split} 
\end{equation}
where $j, l \in \mathbb{Z}$, $j, l \ge 0$, $2j+l \leq m$, $F=\trans(f,g)$ and 
\begin{equation*}
(u_1, u_2)_\gamma=\frac{1}{\gamma^2}(\phi_1, \phi_2)+(w_1, w_2) 
\end{equation*} 
for $u_j=\trans(\phi_j, w_j)$ $(j=1,2);$

\begin{gather}
\begin{split} \label{3}
&\frac{1}{2\gamma^2}\frac{d}{dt}\|T_{j,l}\partial_{x_n}^{k+1}\phi\|_{L^2(\Omega_{per})}^2 +\frac{1}{2(\nu+\tilde{\nu})}\|T_{j,l}\partial_{x_n}^{k+1}\phi\|_{L^2(\Omega_{per})}^2 +\frac{\nu+\tilde{\nu}}{4\gamma^4} \|T_{j,l}\partial_{x_n}^{k+1}\dot{\phi}\|_{L^2(\Omega_{per})}^2  \\
&\leq C\Bigl\{\frac{1}{\gamma^2}(T_{j,l}\partial_{x_n}^{k+1}(w\cdot \nabla \phi),T_{j,l}\partial_{x_n}^{k+1}\phi)+\frac{\nu+\tilde{\nu}}{\gamma^4} \no{T_{j,l}\partial_{x_n}^{k+1}\tilde{f}}{2}{L^2(\Omega_{per})}\\
&\quad +\frac{1}{\nu+\tilde{\nu}}\no{T_{j,l}\partial_{x_n}^{k}g_n}{2}{L^2(\Omega_{per})}
+\frac{1}{\nu+\tilde{\nu}}\no{T_{j+1,l}\partial_{x_n}^{k}w_n}{2}{L^2(\Omega_{per})} \\
& \quad +\frac{\nu^2}{\nu+\tilde{\nu}}\no{T_{j,l+1}\partial_{x_n}^{k}\nabla w}{2}{L^2(\Omega_{per})}\Bigr\},  
\end{split} 
\end{gather}
where $j, l, k \in \mathbb{Z}$, $j, l, k \ge 0$, $2j+l+k \leq m-1$, 
\begin{gather}
\begin{split} \label{4}
&\frac{\nu^2}{\nu+\tilde{\nu}} \|T_{j,l}\partial_{x}^{k+2}w\|_{L^2(\Omega_{per})}^2 +\frac{1}{\nu+\tilde{\nu}}\|T_{j,l}\partial_{x}^{k+1}\phi\|_{L^2(\Omega_{per})}^2  \\
& \quad \leq C\Bigl\{\frac{\nu+\tilde{\nu}}{\gamma^4}\|T_{j,l}\tilde{f}\|_{H^{k+1}(\Omega_{per})}^2+\frac{\nu+\tilde{\nu}}{\gamma^4}\|T_{j,l}\dot{\phi}\|_{H^{k+1}(\Omega_{per})}^2 \\ 
& \quad \quad  +\frac{1}{\nu+\tilde{\nu}}\|T_{j,l}g\|_{H^{k}(\Omega_{per})}^2+\frac{1}{\nu+\tilde{\nu}}\|T_{j,l}\partial_t w\|_{H^{k}(\Omega_{per})}^2 \Bigr\}, 
\end{split} 
\end{gather}
where $j, l, k \in \mathbb{Z}$, $j, l, k \ge 0$, $2j+l+k \leq m-1$, 
\begin{equation} \label{5}
\begin{split}
\frac{\nu+\tilde{\nu}}{\gamma^4}\| \partial_t^{j+1}\phi\|^2_{L^2(\Omega_{per})} 
 \leq C\Bigl\{&(\nu+\tilde{\nu})\|\partial_t^j\mathrm{div} w\|^2_2+\frac{\nu+\tilde{\nu}}{\gamma^4}\|\partial_t^j\tilde{f}\|^2_{L^2(\Omega_{per})} \\ 
& +\frac{\nu+\tilde{\nu}}{\gamma^4}(\partial_t^j(w\cdot\nabla\phi),\partial_t^{j+1}\phi)\Bigr\}, 
\end{split}
\end{equation}
where $j=0,1$.
\end{prop} 

Combining the basic estimates given in Proposition \ref{basicest}, we have the following $H^m$-energy estimate ($m=2,4$) in a similar argument to that in \cite[Section 5]{E-K}. 

\begin{prop} \label{energyest}
There exists positive constant $\nu_0$ such that if $\frac{\nu^2}{\nu+\tilde{\nu}}\geq\nu_0$, then 
\begin{gather}
\frac{d}{dt} \tilde{E}_m(u)+D_m(u)\leq N_m(u) \label{energy}
\end{gather}
for $m=2,4$.
Here 
\begin{gather*}
C^{-1}\left(\frac{1}{\gamma^2}[\![ \phi ]\!]^2_m +\sum_{2j+l\leq m}\|T_{j,l} w\|_{L^2(\Omega_{per})}^2\right) 
\leq \tilde{E}_m(u) \leq C\left(\frac{1}{\gamma^2}[\![ \phi ]\!]^2_m+\sum_{2j+l\leq m}\|T_{j,l} w\|_{L^2(\Omega_{per})}^2\right).
\end{gather*}
and 
\begin{gather*}
\begin{split}
N_m(u) = &C\Big\{\frac{\nu+\tilde{\nu}}{\gamma^4}[\![ \tilde{f} ]\!]^2_m+\frac{1}{\nu+\tilde{\nu}}[\![ g ]\!]^2_{m-1}+\frac{1}{\gamma^2}\sum_{2j+k\leq m}(\partial_t^j\partial_x^k(w\cdot\nabla\phi), \partial_t^j\partial_x^k\phi) \\
&+\frac{1}{\gamma^2}\sum_{2j+l\leq m}(T_{j,l}\tilde{f},T_{j,l}\phi)+\sum_{2j+l\leq m}(T_{j,l}g,T_{j,l}w) + \frac{\nu+\tilde{\nu}}{\gamma^4}(\partial^{\frac{m}{2}-1}_t(w\cdot\nabla\phi),\partial_t^{\frac{m}{2}}\phi)\Big\},
\end{split}
\end{gather*}
where $C$ is independent of $\nu$, $\tilde{\nu}$, $\gamma$ and $S$.
\end{prop}

\hspace{-6mm}{\bf Proof.} \ 
We prove (\ref{energy}) in the case of $m=4$. 
Let $b_j$ ($j=1,\cdots, 9$) be positive numbers independent of $\nu$, $\tilde{\nu}$ and $\gamma$ and consider the following equality:
\begin{gather*} 
\displaystyle 
\sum_{2j+l\leq 4}(\ref{1})
+\sum_{2j+l\leq 3}\left\{b_1\times(\ref{3})_{k=0}+b_2\times(\ref{4})_{k=0}\right\} \\
+\sum_{2j+l\leq 2}\left\{b_3\times(\ref{3})_{k=1}+b_4\times(\ref{4})_{k=1}\right\}
+\sum_{2j+l\leq 3}\left\{b_5\times(\ref{3})_{k=2}+b_6\times(\ref{4})_{k=2}\right\} \\
+b_7\times(\ref{3})_{k=3}+b_8\times(\ref{4})_{k=3} +b_9\times (\ref{5}). 
\end{gather*}
As in \cite[Section 5]{E-K}, taking $b_j$ ($j=1,\cdots,9$) suitably small, if $\nu\geq1$, we can obtain 
\begin{equation}
\frac{d}{dt}\tilde{E}_4(u)+D_4(u)\leq N_4(u), \label{E1}
\end{equation}
where there exists a constant $C$ such that 
\begin{gather*}
C^{-1}\left(\frac{1}{\gamma^2}[\![ \phi ]\!]^2_4+\sum_{2j+l\leq4}\|T_{j,l} w\|_{L^2(\Omega_{per})}^2\right) 
\leq \tilde{E}_4(u) \leq C\left(\frac{1}{\gamma^2}[\![ \phi ]\!]^2_4+\sum_{2j+l\leq4}\|T_{j,l} w\|_{L^2(\Omega_{per})}^2\right).
\end{gather*}
This proves (\ref{energy}) for $m=4$. The case $m = 2$ can be proved in a similar manner.
This completes the proof. \qed

\vspace{2ex}

We next estimate  $N_4(u)$.

\begin{prop}\label{nonlinear}
There exist positive constants $\nu_0$, $\gamma_0$, $a$, $C_0$ and $C_1$ such that 
if $\frac{\nu^2}{\nu+\tilde{\nu}}\geq \nu_0$, $\frac{\gamma^2}{\nu+\tilde{\nu}}\geq \gamma_0$ and $E_4(u)\leq 1$, then the following estimate holds.
\begin{align*}
N_4(u)\leq \frac{1}{4} D_4(u)+C_0\frac{1}{\nu}S^2\sum_{j=0}^2\|\partial_t^jG\|^2_{H^{3-2j}(\Omega_{per})}+C_1\frac{\gamma}{\nu}\sqrt{E_4(u)}D_4(u). 
\end{align*}
\end{prop}

Proposition \ref{nonlinear} is an immediate consequence of the following Proposition \ref{estnorm}.
\begin{prop}\label{estnorm}
Let $j$ and $k$ be non-negative integer with $2j+k\leq4$.
If $\frac{\nu^2}{\nu+\tilde{\nu}}\geq \nu_0$, $\frac{\gamma^2}{\nu+\tilde{\nu}}\geq \gamma_0$ and $E_4(u)\leq 1$, then the following estimates hold.
\begin{enumerate}
\item[{\rm (i)}]  $\displaystyle \frac{\nu+\tilde{\nu}}{\gamma^4}[\![\tilde{f}]\!]_4^2 
\leq \frac{C}{\gamma^2}E_4(u)D_4(u)$ 
\item[{\rm (ii)}] $\displaystyle \frac{\nu+\tilde{\nu}}{\gamma^4}\|\partial_t(w\cdot\nabla\phi)\|_{L^2(\Omega_{per})}^2 
\leq C\frac{\left(\nu+\tilde{\nu}\right)^2}{\gamma^4}E_4(u)D_4(u)$ 
\item[{\rm (iii)}] $\displaystyle \frac{1}{\nu+\tilde{\nu}}[\![g]\!]_3^2 \leq 
\frac{1}{\nu+\tilde{\nu}}S^2[\![G]\!]_3^2
+C\left(\frac{1}{\nu^2}+\frac{1}{\gamma^2} \right)E_4(u)D_4(u)$
\item[{\rm (iv)}] $\displaystyle \frac{1}{\gamma^2}|(\partial_t^j\partial_x^k(w\cdot\nabla\phi),\partial_t^j\partial_x^k\phi)| 
\leq C\frac{\gamma}{\nu} \sqrt{E_4(u)}D_4(u)$
\item[{\rm (v)}] $\displaystyle \frac{1}{\gamma^2}|(\partial_t^j\partial_{x'}^k\tilde{f},\partial_t^j\partial_{x'}^k\phi)|
\leq C\frac{\gamma}{\nu} \sqrt{E_4(u)}D_4(u)$
\item[{\rm (vi)}] $\displaystyle |(\partial_t^j\partial_{x'}^kg, \partial_t^j\partial_{x'}^kw)| 
\leq \frac{1}{4} D_4(u)+\frac{C}{\nu}S^2\sum_{j=0}^2\|\partial_t^jG\|^2_{H^{3-2j}(\Omega_{per})}+C\left(\frac{1}{\nu}+\frac{1}{\gamma}\right)\sqrt{E_4(u)}D_4(u)$
\end{enumerate}
Here $C$ is some positive constant independent of $\nu$, $\tilde{\nu}$, $\gamma$ and $S$.
\end{prop}

To prove Proposition \ref{estnorm}, we use the following Sobolev inequalities. 
\begin{lem}\label{sobolev} 
Let $2\leq p\leq 6$. Then the inequality 
\begin{align*}
\| f\|_{L^p(\Omega_{per})}\leq C\| f\|_{H^1(\Omega_{per})}
\end{align*}
holds for $f\in H^1(\Omega_{per})$. 
Furthermore, the inequality  
\begin{gather*}
\| f\|_{L^\infty(\Omega_{per})}\leq C\|f\|_{H^2(\Omega_{per})}
\end{gather*}
holds for $f\in H^2(\Omega_{per})$.
\end{lem}

\hspace{-6mm}{\bf Proof of Proposition \ref{estnorm}.} 
By straightforward computations based on Lemma \ref{sobolev}, we have 
\begin{gather} 
\frac{\nu+\tilde{\nu}}{\gamma^4}[\![\tilde{f}]\!]_4^2 \leq C\frac{1}{\gamma^2}\left(\frac{\nu+\tilde{\nu}}{\nu}\right)^2E_4(u)D_4(u), \label{nonlinearest1}
\end{gather}
Using (\ref{assump mu}) we obtain the desired estimate (i) from (\ref{nonlinearest1}). 
The estimate (ii) is a direct consequence of Lemma \ref{sobolev}. 

As for the estimate (iii), we first make the following observation. 
We set 
\begin{equation*}
F\left(\frac{\phi}{\gamma^2}\right)=
\begin{cases}\displaystyle 
\frac{\phi}{\gamma^2+\phi}, \\ \displaystyle 
-\frac{1}{\gamma^2+\phi}\left(\frac{1}{\gamma^2}p^{(1)'}(\gamma^{-2}\phi)\phi^2+2p^{(1)}(\gamma^{-2}\phi)\phi-\phi\right).
\end{cases}
\end{equation*}
We see from Lemma \ref{sobolev} that there exist $\nu_{0}$ and $\gamma_{0}$ such that if $E_{4}(u) \le 1$, then $|\phi|\leq \frac{\gamma^2}{2}$ and  
\begin{gather*}\label{estF}
\left\|F\left(\frac{\phi}{\gamma^2}\right)\right\|_{L^\infty(\Omega_{per})}\leq \frac{C}{\gamma^2}\|\phi\|_{L^\infty(\Omega_{per})}, \quad
\left\|F^{(k)}\left(\frac{\phi}{\gamma^2}\right)\right\|_{L^\infty(\Omega_{per})}\leq C
\end{gather*}
for $k=1,2,3$.
By these estimates, together with Lemma \ref{sobolev}, we can obtain the following estimate:  
\begin{gather}
\frac{1}{\nu+\tilde{\nu}}[\![g]\!]_3^2 \leq \frac{1}{\nu+\tilde{\nu}}S^2[\![G]\!]^2_3+C\left(\frac{1}{\nu^2}+\frac{1}{\gamma^2}\left(\frac{\nu+\tilde{\nu}}{\nu}\right)\right)E_4(u)D_4(u). \label{nonlinearest3} 
\end{gather}
Using (\ref{assump mu}) we obtain the desired estimate (iii) from (\ref{nonlinearest3}). 

As for (iv), we have 
\begin{align*}
|(\partial_t^j\partial_x^k(w\cdot\nabla\phi),\partial_t^j\partial_x^k\phi)| 
&= |((w\cdot\nabla\partial_t^j\partial_x^k\phi),\partial_t^j\partial_x^k\phi)
+([\partial_t^j\partial_x^k,w\cdot\nabla]\phi,\partial_t^j\partial_x^k\phi)| \\
&= |-(\mathrm{div}w ,\frac{1}{2}|\partial_t^j\partial_x^k\phi|^2)
+([\partial_t^j\partial_x^k,w\cdot\nabla]\phi,\partial_t^j\partial_x^k\phi)|. 
\end{align*}
Here $[A,B]$ denotes the commutator of $A$ and $B$: $[A,B]f=A(Bf)-B(Af)$. 
Applying Lemma \ref{sobolev} we obtain the estimate (iv). 
A direct application of Lemma \ref{sobolev} gives the estimate (v). 
As for (vi), we also have 
\begin{gather}
\begin{split}
|(\partial_t^j\partial_{x'}^k g,\partial_t^j\partial_{x'}^kw)| 
\leq &\frac{1}{4} D_4(u)+C\frac{\nu+\tilde{\nu}}{\nu^2}S^2\sum_{j=0}^2\|\partial_t^jG\|^2_{H^{3-2j}(\Omega_{per})} \\
&+C\left\{\frac{\nu+\tilde{\nu}}{\nu^2}+\frac{1}{\gamma}\left(\frac{\nu+\tilde{\nu}}{\nu}\right)^2\right\}\sqrt{E_4(u)}D_4(u). 
\end{split}\label{nonlinearest6}
\end{gather}
Using (\ref{assump mu}) we obtain the desired estimate (vi) from (\ref{nonlinearest6}).  \qed 

\begin{rem}\label{upperbound}
The main reason we make the assumption $(\ref{assump mu})$ lies in Proposition $\ref{estnorm}$. 
In its proof, we use $(\ref{assump mu})$ to obtain the estimates $({\rm i})$, $({\rm iii})$ and $({\rm v})$ of Proposition $\ref{estnorm}$ from $(\ref{nonlinearest1})$, $(\ref{nonlinearest3})$ and $(\ref{nonlinearest6})$, respectively. 
One could prove the existence of the time periodic solution without the assumption $(\ref{assump mu})$, by using  Proposition $\ref{estnorm}$ with $({\rm i})$, $({\rm iii})$ and $({\rm v})$ replaced by $(\ref{nonlinearest1})$, $(\ref{nonlinearest3})$ and $(\ref{nonlinearest6})$, in which case, however, the condition on $S$ would become more complicated. 
We thus consider the problem under the assumption $(\ref{assump mu})$. 
\end{rem}

\vspace{2ex}
We prepare the following lemma which is proved by using extension operator. 

\begin{lem}\label{extension}
Let $k \in \mathbb{N}$. There exists a positive constant $C_E$ such that 
\begin{align*}
\|f(t)\|_{H^{k}(\Omega_{per}) }^{2} & \le C_E \left\{ e^{ - b(t - s)}\|f(s)\|_{H^{k}(\Omega_{per})}^{2} + b \int_{s}^{t}e^{-b(t -\tau)} \|f(s)\|_{H^{k}(\Omega_{per})}^{2}\,d\tau \right. \\[1ex]
& \quad \quad \left. + 2 \int_{s}^{t}e^{-b(t -\tau)} \|\partial_{t} f(s)\|_{H^{k-1}(\Omega_{per})}\| f(s)\|_{H^{k+1}(\Omega_{per})}\,d\tau \right\}
\end{align*}
for any $0\le s \le t$, $b \ge $ and $f\in L^{2}_{loc}([0, \infty); H^{k+1}(\Omega_{per}) )\cap H^{1}_{loc}([0, \infty); H^{k-1}(\Omega_{per}))$. 
\end{lem}

We are now in a position to prove Proposition \ref{apriori}.

\vspace{2ex}
\noindent 
{\bf Proof of Proposition \ref{apriori}.} 
Assume that 
\begin{equation}\label{est E4}
E_4(u(t))\leq 6\frac{C_0(1+C_E)}{\nu}S^2\frac{1}{1-e^{-a\frac{\nu+\tilde{\nu}}{\gamma^2}}}
\end{equation} 
for $t\in [0,T]$ and $S$ satisfies 
\begin{equation}
S\leq\frac{1}{10 \sqrt{C_0(1 + C_E)}}\frac{\nu}{\sqrt{\nu+\tilde{\nu}}}\sqrt{1-e^{-a\frac{\nu+\tilde{\nu}}{\gamma^2}}}\min\left\{1, \frac{1}{C_1}\frac{\nu}{\gamma} \right\}, \label{S}
\end{equation}
where $C_0$ and $C_1$ are the constants given in Proposition \ref{nonlinear}. 
We then have $E_4(u(t))<1$ for $t\in[0,T]$. 
It then follows from Propositions \ref{energyest} and \ref{nonlinear} that 
\begin{equation*}
\frac{d}{dt} \tilde{E}_4(u)+\frac{3}{4}D_4(u)\leq \frac{C_0}{\nu}S^2\sum_{j=0}^2\|\partial_t^jG\|^2_{H^{3-2j}(\Omega_{per})}+C_1\frac{\gamma}{\nu}\sqrt{E_4(u)}D_4(u).
\end{equation*} 
By the Poincar{\'e} inequality, we have 
\begin{equation*}
D_4(u)\geq 8a\frac{\nu+\tilde{\nu}}{\gamma^2}E_4(u), 
\end{equation*}
and hence, 
\begin{equation*}
\frac{d}{dt} \tilde{E}(u)+2a\frac{\nu+\tilde{\nu}}{\gamma^2}E_4(u)+\frac{1}{2}D_4(u)\leq \frac{C_0}{\nu}S^2\sum_{j=0}^2\|\partial_t^jG\|^2_{H^{3-2j}(\Omega_{per})}+C_1\frac{\gamma}{\nu}\sqrt{E_4(u)}D_4(u).
\end{equation*}
Using (\ref{est E4}), (\ref{S}) and the relation $\tilde{E}_{4}(u) \le E_{4}(u)$, we find that 
\[
\frac{d}{dt}\tilde{E}_4(u)+a\frac{\nu+\tilde{\nu}}{\gamma^2}(\tilde{E}_{4}(u) + E_4(u))+\frac{1}{4}D_4(u)\leq \frac{C_0}{\nu}S^2\sum_{j=0}^2\|\partial_t^jG\|^2_{H^{3-2j}(\Omega_{per})}. 
\]
This gives 
\begin{equation}\label{est global sol}
\begin{split}
&\tilde{E}_{4}(u(t))+ a\frac{\nu + \tilde{\nu}}{\gamma^{2}} \int_{0}^{t} e^{-a\frac{\nu + \tilde{\nu}}{\gamma^{2}}(t - s)} E_{4}(u) \, ds + \frac{1}{4}\int_0^t e^{-a\frac{\nu+\tilde{\nu}}{\gamma^2}(t-s)}D_4(u(s))\,ds  \\
& \leq e^{-a\frac{\nu+\tilde{\nu}}{\gamma^2}t} \tilde{E}_4(u_0) 
+\frac{C_0}{\nu}S^2\int^t_0e^{-a\frac{\nu+\tilde{\nu}}{\gamma^2}(t-s)}\sum_{j=0}^2\|\partial_t^jG(s)\|^2_{H^{3-2j}(\Omega_{per})}\,ds. 
\end{split}
\end{equation}
By Lemma \ref{extension} with $b=a\frac{\nu + \tilde{\nu}}{\gamma^{2}}$, we see, by taking $\nu_{0}$ suitably large, that  
\begin{align*}
\nom{\partial_{x_n}^{2}w(t)}{2}{2} 
& \le C_E\left\{ e^{-a\frac{\nu+\tilde{\nu}}{\gamma^{2}}t} \nom{\partial_{x_n}^{2}w_{0}}{2}{2} + a\frac{\nu + \tilde{\nu}}{\gamma^{2}} \int_{0}^{t} e^{-a\frac{\nu + \tilde{\nu}}{\gamma^{2}}(t - s)} E_{4}(u) \, ds \right. \\[1ex]
& \quad \quad \left. + \frac{1}{4} \int_0^t e^{-a\frac{\nu+\tilde{\nu}}{\gamma^2}(t-s)}D_4(u(s))\,ds \right\}. 
\end{align*}
By adding this to $(1+C_E) \times (\ref{est global sol})$ and using Lemma \ref{intlemma}, we deduce that  
\begin{align*}
& E_4(u(t))+ \frac{1}{4} \int_0^te^{-a\frac{\nu+\tilde{\nu}}{\gamma^2}(t-s)}D_4(u(s))\,ds  \\[1ex] 
& \quad \le (1+C_E) \left\{e^{-a\frac{\nu+\tilde{\nu}}{\gamma^2}t}E_4(u_0) 
+2\frac{C_0}{\nu}S^2\frac{1}{1-e^{-a\frac{\nu+\tilde{\nu}}{\gamma^2}}} \right\}
\end{align*}
for all $t\in[0,T]$. 
This completes the proof. \qed

\subsection{Existence of time periodic solution} \label{existtimeper}

We first consider the $H^2$-energy estimate for the difference of the solution of (\ref{prob3})-(\ref{IC}). 

Let $u_{0j}$ $(j=1,2)$ and $G$ satisfy \eqref{cc0} and the compatibility conditions (\ref{cc1}) and (\ref{cc2}). 
Assume that $u_{0j}$ satisfy 
\begin{equation}
E_4(u_{0j})\leq 4\frac{C_0(1+C_E)}{\nu}S^2\frac{1}{1-e^{-a\frac{\nu+\tilde{\nu}}{\gamma^2}}} \quad (j=1,2). \label{Eu0j}
\end{equation}
Let $u_j$ ($j=1,2$) be the solutions of  (\ref{prob3})-(\ref{IC}) with $u_0=u_{0j}$ obtained by Proposition \ref{est} with $u$ and $u_0$ replaced by $u_j$ and $u_{0j}$ ($j=1,2$), respectively. 

We set $\tilde{u}=u_1-u_2=\trans(\phi_1-\phi_2,w_1-w_2)$. 
Then $\tilde{u}=\trans(\tilde{\phi}, \tilde{w})$ satisfies 
\begin{equation} \label{eq2-1}
\partial_t \tilde{\phi} + \gamma^2 \mathrm{div} \tilde{w} = f, 
\end{equation}
\begin{equation}\label{eq2-2}
\partial_t \tilde{w} -\nu \Delta \tilde{w}- \tilde{\nu}\nabla \mathrm{div} \tilde{w} + \nabla \tilde{\phi} = g,\\  
\end{equation}
where 
\begin{gather*}
f=-w_1\cdot\nabla\tilde{\phi}-\tilde{f}, \\
\begin{split}
g=-&\{ w_1\cdot\nabla \tilde{w}+\tilde{w}\cdot\nabla w_2+\frac{\phi_1}{\gamma^2+\phi_1}\Delta \tilde{w}+\frac{1}{\gamma^2+\phi_1}\frac{\gamma^2}{\gamma^2+\phi_2}\Delta w_2 \tilde{\phi} \\
&+F\left(\frac{\phi_1}{\gamma^2}\right)\nabla\tilde{\phi}+F^{(1)}\left(\frac{\phi_1}{\gamma^2},\frac{\phi_2}{\gamma^2}\right)\nabla\phi_2\frac{1}{\gamma^2}\tilde{\phi}\}.
\end{split}
\end{gather*}
Here 
\begin{gather*}
\tilde{f}=\tilde{w}\cdot\nabla\phi_2+\phi_1\mathrm{div}\tilde{w}+\phi\mathrm{div}w_2, \\
F\left(\frac{\tilde{\phi}}{\gamma^2}\right)=-\frac{1}{\gamma^2+\tilde{\phi}}\left(\frac{1}{\gamma^2}p^{(1)'}\left(\frac{\tilde{\phi}}{\gamma^2}\right)\tilde{\phi}^2+2p^{(1)}\left(\frac{\tilde{\phi}}{\gamma^2}\right)\tilde{\phi}-\tilde{\phi}\right), \\
F^{(1)}\left(\frac{\phi_1}{\gamma^2}, \frac{\phi_2}{\gamma^2}\right)=\int^1_0F'\left(\theta\frac{\phi_1}{\gamma^2}+(1-\theta)\frac{\phi_2}{\gamma^2} \right) d\theta. 
\end{gather*}

Similarly to the previous section one can obtain 
\begin{equation}
\frac{d}{dt}\tilde{E}_2(\tilde{u})+D_2(\tilde{u})\leq N_2(\tilde{u}). \label{tildeenergy}
\end{equation}

One can also obtain the following estimate for $N_2(\tilde{u})$ in a similar manner to the previous section.
\begin{prop}\label{nonlinear2}
If $\frac{\nu^2}{\nu+\tilde{\nu}}\geq \nu_0$, $\frac{\gamma^2}{\nu+\tilde{\nu}}\geq \gamma_0$ and $S$ satisfies $(\ref{S})$, then there exists positive constant $C_2$ independent of $\nu$, $\tilde{\nu}$, $\gamma$ and $S$ such that 
\begin{equation*}
N_2(\tilde{u})\leq C_2\frac{\gamma^2}{\nu+\tilde{\nu}}\left(\sqrt{E_4(u_1)}+\sqrt{E_4(u_2)}\right)D_2(\tilde{u}).
\end{equation*}
\end{prop}

In fact, Proposition \ref{nonlinear2} is an immediate consequence of the following Proposition \ref{estinner2}. 
Here we also use the assumption (\ref{assump mu}) as in the proof of Proposition \ref{estnorm}.
\begin{prop}\label{estinner2}
If $\frac{\nu^2}{\nu+\tilde{\nu}}\geq \nu_0$, $\frac{\gamma^2}{\nu+\tilde{\nu}}\geq \gamma_0$ and $S$ satisfies $(\ref{S})$, then 
\begin{enumerate}
\item[${\rm (i)}$] $\displaystyle \frac{\nu+\tilde{\nu}}{\gamma^4}[\![\tilde{f}]\!]_2^2 
\leq C
\left(E_4(u_1)+E_4(u_2)\right)D_2(\tilde{u})$, 
\item[${\rm (ii)}$] $\displaystyle \frac{1}{\nu+\tilde{\nu}}[\![g]\!]_1^2
\leq C
\left(E_4(u_1)+E_4(u_2)\right)D_2(\tilde{u})$, 
\item[${\rm (iii)}$] $\displaystyle \frac{\nu+\tilde{\nu}}{\gamma^4}\|w_1\cdot\nabla\tilde{\phi}\|^2_{L^2(\Omega_{per})} 
\leq \frac{(\nu+\tilde{\nu})^2}{\gamma^4}
E_4(u_1)D_2(\tilde{u})$, 
\item[${\rm (iv)}$] $\displaystyle \frac{1}{\gamma^2}|(\partial_t^j\partial_x^k(w_1\cdot\nabla\tilde{\phi}),\partial_t^j\partial_x^k\tilde{\phi})| 
\leq C\frac{\gamma^2}{\nu+\tilde{\nu}}
\left(\sqrt{E_4(u_1)}+\sqrt{E_4(u_2)}\right)D_2(\tilde{u})$, 
\item[${\rm (v)}$] $\displaystyle \frac{1}{\gamma^2}|(\partial_t^j\partial_x^k\tilde{f},\partial_t^j\partial_x^k\tilde{\phi})| 
\leq C\frac{\gamma^2}{\nu+\tilde{\nu}}
\left(\sqrt{E_4(u_1)}+\sqrt{E_4(u_2)}\right)D_2(\tilde{u})$, 
\item[${\rm (vi)}$] $\displaystyle  |(\partial_t^j\partial_x^kg,\partial_t^j\partial_x^k\tilde{w})| 
\leq C\left(1+\frac{\gamma}{\nu} \right) 
\left(\sqrt{E_4(u_1)}+\sqrt{E_4(u_2)}\right)D_2(\tilde{u})$, 
\end{enumerate}
where $C$ is some constant independent of $\nu$, $\tilde{\nu}$, $\gamma$ and $S$. 
\end{prop}

Proposition \ref{estinner2} can be proved in a similar manner to Proposition \ref{estnorm}. 

We now establish the $H^2$-energy estimate for $\tilde{u}$.
\begin{prop} \label{tildeest}
Let $u_j$ $(j=1,2)$ be the solutions of $(\ref{prob3})$-$(\ref{IC})$ with $u_0=u_{0j}$ $(j=1,2)$.  
If $\frac{\nu^2}{\nu+\tilde{\nu}}\geq \nu_0$, $\frac{\gamma^2}{\nu+\tilde{\nu}}\geq \gamma_0$, then the following estimate holds. 
There exists a positive constant $\varepsilon_1$ such that if $S\leq \varepsilon_{1}\frac{\nu^2}{\gamma^2\sqrt{\nu+\tilde{\nu}}}\sqrt{1-e^{-a\frac{\nu+\tilde{\nu}}{\gamma^2}}}$, then  
\begin{equation}
E_2(\tilde{u}(t))+\frac{1}{4}\int^t_0e^{-a\frac{\nu+\tilde{\nu}}{\gamma^2}(t-s)}D_2(\tilde{u}(s))\,ds \leq (1+C_E) e^{-a\frac{\nu+\tilde{\nu}}{\gamma^2}t}E_2(u_{01}-u_{02}). \label{Etilde}
\end{equation}
\end{prop}

\noindent 
{\bf Proof.} By (\ref{tildeenergy}) and Proposition \ref{nonlinear2}, we have   
\begin{equation*}
\frac{d}{dt} \tilde{E}_2(\tilde{u})+2a \frac{\nu+\tilde{\nu}}{\gamma^2}E_2(\tilde{u})+D_2(\tilde{u})
\leq C_2\frac{\gamma^2}{\nu+\tilde{\nu}}\left( \sqrt{E_4(u_1)}+\sqrt{E_4(u_2)}\right)D_2(\tilde{u}).
\end{equation*}

It then follows that if $\frac{\nu^2}{\nu+\tilde{\nu}}\geq \nu_0$, $\frac{\gamma^2}{\nu+\tilde{\nu}}\geq \gamma_0$ and $S$ satisfies 
\begin{equation}\label{S1}
S\leq \frac{1}{10 C_2\sqrt{C_0(1 + C_E)}}\frac{\nu^2}{\gamma^2\sqrt{\nu+\tilde{\nu}}}\sqrt{1-e^{-a\frac{\nu+\tilde{\nu}}{\gamma^2}}}, 
\end{equation}
then 
\begin{equation*}
\frac{d}{dt}\tilde{E}_2(\tilde{u})+2a\frac{\nu+\tilde{\nu}}{\gamma^2}E_2(\tilde{u})+\frac{1}{4}D_2(\tilde{u})\leq 0.
\end{equation*}
Proposition \ref{tildeest} now follows from this inequality and Lemma \ref{extension} in a similar manner to the proof of Proposition \ref{apriori}. 
This completes the proof.  \qed

\vspace{2ex}

We now show the existence of a time-periodic solution of (\ref{prob3})-(\ref{meanvalue}). 

\vspace{2ex}

\hspace{-6mm}{\bf Proof of Proposition \ref{prop1}.}
Let $\phi^\flat_0=0$ and $w_0^\flat\in H^4(\Omega_{per})\cap H^1_0(\Omega_{per})$ be the solution to 
\begin{gather*} 
-\nu \Delta w- \tilde{\nu}\nabla \mathrm{div} w = g(0, w, G(0)). 
\end{gather*}
The existence of $w_{0}^{\flat}$ follows from the standard elliptic theory. 
Furthermore, one can prove that $u^{\flat}=\trans(\phi_{0}^{\flat}, w_{0}^{\flat})$ satisfies 
\begin{gather}\label{eststationary}
E_4(u^\flat_0)\leq \frac{2C_0 (1+ C_E) }{\nu}S^2\frac{1}{1-e^{-a\frac{\nu+\tilde{\nu}}{\gamma^2}}}. 
\end{gather}
This can be seen by a similar energy method as that in Section 5.1 without time derivatives and by applying Lemma \ref{extension} to $G$.
By Proposition \ref{est}, we have the global solutions $u^\flat(t)$ of (\ref{prob3})-(\ref{IC}) with $u_0=u^\flat_0$ and $u^\flat(t)$ satisfies 
\begin{equation*}
E_4(u^\flat(t))\leq \frac{4C_0(1+C_E)}{\nu}S^2\frac{1}{1-e^{-a\frac{\nu+\tilde{\nu}}{\gamma^2}}} \quad (t\geq0). 
\end{equation*}

We next consider the functions $u_1$ and $u_2$ defined by 
\begin{equation*}
u_1(t)=u^\flat(t), \quad u_2(t)=u^\flat(t+(m-n)),
\end{equation*}
where $m,n \in\mathbb{N}$ with $m>n$.
As in the proof of Proposition \ref{L2 time}, we can show that 
\begin{equation*}
E_2(u^\flat(n)-u^\flat(m))\rightarrow 0 \quad (n\rightarrow \infty), 
\end{equation*} 
and  it holds that there exists $\tilde{u}^\flat \in H_*^4(\Omega_{per})\times (H^4(\Omega_{per})\cap H^1_0(\Omega_{per}))$ such that 
$u^\flat(m)$ converges to $\tilde{u}^\flat$ strongly in $H_*^2(\Omega_{per})\times(H^2(\Omega_{per})\cap H^1_0(\Omega_{per}))$ and weakly in $H_*^4(\Omega_{per})\times(H^4(\Omega_{per})\cap H^1_0(\Omega_{per}))$, and $\tilde{u}^\flat$ satisfies   
\begin{gather*}
E_4(\tilde{u}^\flat)\leq 4\frac{C_0(1+C_E)}{\nu}S^2\frac{1}{1-e^{-a\frac{\nu+\tilde{\nu}}{\gamma^2}}}.
\end{gather*}
We therefore see from Proposition \ref{est} that there exists a unique global solution $u\in\cap_{j=0}^2C^j([0,\infty);H^{4-2j}(\Omega_{per})\times H^{4-2j}(\Omega_{per}))\cap H^j([0,\infty); H^{4-2j}(\Omega_{per})\times H^{5-2j}(\Omega_{per}))$ of (\ref{prob3})-(\ref{IC}) in with $u_0=\tilde{u}^\flat$. It then follows from the argument by Valli \cite{Valli} that $u$ is a time-periodic solution of (\ref{prob3})-(\ref{meanvalue}) satisfying 
\begin{gather*}
E_4(u(t))\leq \frac{C}{\nu}S^2\frac{1}{1-e^{-a\frac{\nu+\tilde{\nu}}{\gamma^2}}}, \\
\int_0^1D_4(u(s))\,ds\leq \frac{C}{\nu}S^2\frac{1}{1-e^{-a\frac{\nu+\tilde{\nu}}{\gamma^2}}}.
\end{gather*}
Using the condition (\ref{S1}), we have 
\begin{gather*}
\frac{1}{\gamma^2}\nom{\phi}{2}{4}+\nom{w}{2}{4}\leq C\frac{\nu^2}{\gamma^4}, \\
\int_0^1\left\{\frac{\nu^2}{\nu+\tilde{\nu}}\nom{w(s)}{2}{4}+\frac{1}{\nu+\tilde{\nu}}\nom{\nabla\phi(s)}{2}{3}+\frac{\nu+\tilde{\nu}}{\gamma^4}\|\partial_t^2\phi(s)\|^2_{L^2(\Omega_{per})}\right\}\,ds\leq C\frac{\nu^2}{\gamma^4}.
\end{gather*}
This completes the proof. \qed

\vspace{3ex}
\noindent
{\bf Acknowledgements.} 
The authors would like to thank the reviewers for their valuable comments.  
Y. Kagei was partly supported by JSPS KAKENHI Grant Numbers 16H03947, 16H06339 and 20H00118.
S. Enomoto was partly supported by JSPS KAKENHI Grant Numbers 18J01068.

\end{document}